\documentclass[reqno,12pt,oneside,a4paper]{amsart} 

\usepackage{bbm,enumerate,amsmath,amstext,mathrsfs,a4}
\usepackage[english]{babel}
\usepackage[latin1]{inputenc}

\newcommand{\R}{\mathbbm{R}}
\newcommand{\N}{\mathbbm{N}}
\newcommand{\Q}{\mathbbm{Q}}
\newcommand{\Z}{\mathbbm{Z}}
\newcommand{\C}{\mathbbm{C}}
\newcommand{\CP}{\C P}
\newcommand{\RP}{\R P}
\newcommand{\CPl}{\CP^\infty_{-1}}
\newcommand{\CPll}{\Omega^\infty\Sigma\CPl}
\newcommand{\CPlll}{\Omega^\infty\CPl}
\newcommand{\Th}{\mathbbm{T}\mathrm{h}}

\newcommand{\Diff}{\mathrm{Diff}}
\newcommand{\F}{\mathbbm{F}}

\newcommand{\Tor}{\mathrm{Tor}}
\newcommand{\Cotor}{\mathrm{Cotor}}
\newcommand{\modd}{\backslash\!\!\backslash}

\newcommand{\Ker}{\mathrm{Ker}}
\newcommand{\Cok}{\mathrm{Cok}}
\newcommand{\IM}{\mathrm{Im}}
\newcommand{\moddd}{/\!\!/}

\newcommand{\trf}{\partial}
\renewcommand{\phi}{\varphi}
\renewcommand{\epsilon}{\varepsilon}

\newcommand{\RR}{\mathscr{R}}
\newcommand{\DD}{\mathscr{D}}
\newcommand{\JJ}{\mathscr{J}}
\newcommand{\CC}{\mathscr{C}}
\renewcommand{\AA}{\mathscr{A}}

\newcommand{\Prim}{P}
\newcommand{\Ind}{Q}

\newcommand{\cotensor}{\Box}
\newcommand{\tensor}{\otimes}
\newcommand{\isom}{\cong}

\newcommand{\colim}{\operatorname{colim}}
\newcommand{\holim}{\operatorname{holim}}

\usepackage{amsthm}
\theoremstyle{plain}
\newtheorem{thm}{Theorem}[section]
\newtheorem{prop}[thm]{Proposition}
\newtheorem{lem}[thm]{Lemma}

\newtheorem{cor}[thm]{Corollary}

\theoremstyle{definition}

\newtheorem{defn}[thm]{Definition}

\theoremstyle{remark}

\newtheorem{rem}[thm]{Remark}

\numberwithin{equation}{section}

\usepackage[all]{xy}
\CompileMatrices

\begin{document}

\title{Mod $p$ homology of the stable mapping class group}
\author{Søren Galatius}
\address{Aarhus University, Aarhus, Denmark}
\email{galatius@imf.au.dk}

\begin{abstract}
  We calculate the homology $H_*(\Gamma_{g,n};\F_p)$ of the mapping
  class group $\Gamma_{g,n}$ in the stable range.  The calculation is
  based on Madsen and Weiss' proof of the ``Generalised Mumford
  Conjecture'':  $\Gamma_{g,n}$ has the same homology as a component of
  the space $\CPlll$ in the stable range.
\end{abstract}

\maketitle

\section{Introduction}
\label{sec:introduction}

Let $F_{g,n}$ be an oriented surface of genus $g$ and with $n$
boundary components.  Let $\Diff(F_{g,n}, \partial)$ be the
topological group of orientation preserving diffeomorphisms of
$F_{g,n}$ fixing pointwise a neighbourhood of the boundary.  The
mapping class group is the group $\Gamma_{g,n} =
\pi_0\Diff(F_{g,n},\partial)$ of components.  There are group maps
\begin{equation*}
  \Gamma_{g,n} \to \Gamma_{g,n-1} \quad\text{and}\quad \Gamma_{g,n}
  \to \Gamma_{g+1,n}
\end{equation*}
induced by gluing a disk, resp.\ a torus with two boundary components,
to one of the boundary components of $F_{g,n}$.  By a theorem of Harer
and Ivanov, these maps induce isomorphisms in $H_*(-;\Z)$ for $*\leq
(g-1)/2$, and thus there is a stable range in which the group homology
$H_*(\Gamma_{g,n};\Z)$ is independent of $g$ and $n$.  In this range
it agrees with $H_*(\Gamma_{\infty};\Z)$ where $\Gamma_{\infty} =
\colim_g \Gamma_{g,1}$ is the stable mapping class group.

\subsection{Madsen-Weiss' theorem}
\label{sec:madsen-weiss-theorem}

The Mumford conjecture predicts that
\begin{equation*}
  H^*(\Gamma_\infty;\Q) \isom \Q[\kappa_1, \kappa_2, \dots]
\end{equation*}
for certain classes $\kappa_i\in H^{2i}(\Gamma_\infty)$.  This was
recently proved by Madsen and Weiss, but their result gives more.  To
state the full result, consider the classifying space
$B\Gamma_{\infty}$.  Its homology is the group homology of
$\Gamma_{\infty}$.  By the Quillen plus-construction we get a simply
connected space $B\Gamma_\infty^+$ and a map $B\Gamma_\infty \to
B\Gamma_\infty^+$ inducing an isomorphism in homology.  The
Madsen-Weiss theorem determines the homotopy type of $\Z\times
B\Gamma_\infty^+$ to be that of $\CPlll$.  The space
$\Omega^\infty\CP^\infty_{-1}$ (to be defined below), can be examined
by methods from stable homotopy theory.  In particular it is easy to
calculate $H^*(\Omega^\infty\CP^\infty_{-1};\Q)$.  This implies the
Mumford conjecture.

In this paper we calculate $H_*(\Omega^\infty\CP^\infty_{-1};\F_p)$
for any prime $p$ and hence by the above $H_*(\Gamma_{g,n};\F_p)$ for
$*\leq(g-1)/2$.

\subsection{Outline and statement of results}
\label{sec:outl-stat-results}

Let $p$ be a prime number, and let $H_* = H_*(-;\F_p)$.

Let $L$ be the canonical complex line bundle over $\CP^\infty$ and let
$\CP^\infty_{-1} = \Th(-L)$ be the Thom spectrum of the
$-2$-dimensional virtual bundle $-L$.  Inclusion of a point in
$\CP^\infty$ induces a map $S^{-2} \to \CP^\infty_{-1}$ and the zero
section of the line bundle $L$ induces a map
\begin{equation*}
  \CP^\infty_{-1} = \Th(-L) \to \Th(-L+L) = \Sigma^\infty\CP^\infty_+
\end{equation*}
These fit together into a cofibration sequence
\begin{equation*}
  S^{-2} \to \CP^\infty_{-1} \to \Sigma^\infty\CP^\infty_+
\end{equation*}
and there is an induced fibration sequence of infinite loop spaces
\begin{equation}\label{eq:24}
  \xymatrix{
    {\Omega^\infty\Sigma\CPl} \ar[r]^{\omega} &
    {Q(\Sigma\CP^\infty_+)} \ar[r]^-\partial & QS^0
  }
\end{equation}
where we write $Q = \Omega^\infty\Sigma^\infty$.  This fibration
sequence is the starting point for the calculation of the mod $p$
homology of $\CPll$ and $\CPlll$.

The mod $p$ homology of $Q\Sigma\CP^\infty_+$ and $QS^0$ is completely
known (\cite{AK}, \cite{DL}, \cite{CLM}), as is the induced map
$\partial_*$ in homology (\cite{MMM}).  The first main result of this
paper is a calculation of the Hopf algebra $H_*(\CPll;\F_p)$.  We need
to introduce the following notation to state the results (see
Section~\ref{sec:recollections} for further details).  All Hopf
algebras will be commutative and cocommutative.  The Hopf algebra
cokernel and kernel of a map $f: A\to B$ of Hopf algebras will be
denoted $B\moddd f$ and $A\modd f$, respectively.  $PA$ is the
vectorspace of primitive elements, and $QA$ is the vectorspace of
indecomposable elements.  For a graded vectorspace $V$, $s^{-1}V$ will
denote the desuspension of $V$: $(s^{-1}V)_{\nu-1} = V_{\nu}$.  We
also need to introduce the following functors from vectorspaces to
algebras.  Let $V$ be a non-negatively graded vectorspace and let
$B\subseteq V_0$ be a basis for the degree zero part of $V$.  Let $V =
V^+ \oplus V^-$ be the splitting of $V$ into even and odd dimensional
parts.  Then $E[V^-]$ is the exterior algebra generated by $V^-$ and
$\F_p[V^+]$ is the polynomial algebra generated by $V^+$.  Furthermore
$\F_p[B,B^{-1}]$ is the algebra of Laurent polynomials in the elements
of $B$ and $\F_p[V^+,B^{-1}] = \F_p[V^+]\tensor_{\F_p[B]}
\F_p[B,B^{-1}]$ is the polynomial algebra generated by $V^+$, with the
elements of $B$ inverted.  The \emph{free commutative algebra}
generated by $V$ is $S[V] = E[V^-]\tensor \F_p[V^+,B^{-1}]$.

The calculations use the theory of \emph{homology operations}.  These
are defined on the mod $p$ homology of infinite loop spaces, and are
natural with respect to infinite loop maps, cf.\ \cite{AK}, \cite{DL},
\cite{CLM}.  The basic operations are
\begin{align*}
  \beta^\epsilon Q^s & : H_n(X) \to H_{n+2s(p-1)-\epsilon}(X), &
  &(p>2)\\
  Q^s & : H_n(X) \to H_{n+s}(X), & & (p=2)
\end{align*}
where $\epsilon\in\{0,1\}$ and $s\in\Z_{\geq \epsilon}$.  Given a
sequence $I=(\epsilon_1, s_1, \dots, \epsilon_k, s_k)$ (with all
$\epsilon_i = 0$ if $p=2$) there is an iterated operation $Q^I =
\beta^{\epsilon_1} Q^{s_1} \dots \beta^{\epsilon_k} Q^{s_k}$.  The mod
$p$ homology of $QS^0$ then has the following form: Let $\iota\in
H_0(QS^0)$ be the image of the non-basepoint in $S^0$.  Then
$H_*(QS^0)$ is the free commutative algebra on the set
\begin{equation}
  \label{eq:7}
  \mathbf{T} = \{Q^I \iota \vert \text{$I$ admissible, $e(I) + b(I) >
  0$}\}
\end{equation}
(see section 3 for the definition of $e(I)$, $b(I)$ and the notion of
admissibility).

As a step towards calculating $H_*(\CPll)$ we determine the Hopf
algebra cokernel of $\trf_*: H_*(Q\Sigma\CP^\infty_+) \to H_*(QS^0)$.
In the following theorem, $Q_0S^0 \subseteq QS^0$ is the basepoint
component.  Then $H_*(QS^0) = H_0(QS^0) \tensor H_*(Q_0S^0)$.
\begin{thm}
  \label{thm-2}
  Let $\mathbf{T}$ be as in~\eqref{eq:7}.  Let $H_*(QS^0)^{(0)}$
  denote the subalgebra of $H_*(QS^0)$ generated by the set
  \begin{align*}
    &\{Q^I\iota\in\mathbf{T} \vert \text{all $\epsilon_i = 0$}\} & &
    (p>2)\\
    &\{Q^I\iota\in\mathbf{T} \vert \text{all $s_i$ even}\} & & (p=2)
  \end{align*}
  Then the composite
  \begin{equation*}
    H_*(QS^0)^{(0)} \to H_*(QS^0) \to H_*(QS^0)\moddd\trf_*
  \end{equation*}
  is an isomorphism.  In particular the Hopf algebra
  $H_*(QS^0)\moddd\trf_*$ is concentrated in degrees $\equiv 0
  \pmod{2(p-1)}$.  Similarly for $H_*(Q_0S^0)\moddd\trf_*$.
  Furthermore the dual algebra $H^*(Q_0S^0)\modd\trf^*$ is a
  polynomial algebra.
\end{thm}
\begin{thm}
  \label{thm-1}
  The sequence
  \begin{equation*}
    \xymatrix{
      {\F_p} \ar[r]& {H_*(\CPll)\modd\omega_*} \ar[r]& {H_*(\CPll)}
      \ar[r]^-{\omega_*} &
      {H_*(Q\Sigma\CP^\infty_+)\modd \partial_*} \ar[r]& {\F_p}
    }
  \end{equation*}
  is an exact sequence of Hopf algebras.  It is split but not
  canonically.  Furthermore there is a canonical isomorphism
  \begin{equation*}
    H_*(\CPll)\modd\omega_* \isom E[s^{-1}P(H_*(QS^0)\moddd\partial_*)]
  \end{equation*}
  In particular, $H_*(\CPll)$ is primitively generated and for $p>2$
  it is free commutative.
\end{thm}
Theorem~\ref{thm-2} is an algebraic consequence of the known structure
of $H_*(Q\Sigma\CP^\infty_+)$ and $H_*(QS^0)$ and the induced map
$\trf_*$ in homology.  The proof of Theorem~\ref{thm-1} uses the
Eilenberg-Moore spectral sequence of the fibration
sequence~\eqref{eq:24}.

Next we calculate $H_*(\CPlll)$.
\begin{thm}
  \label{thm-3}
  For $p = 2$, the sequence
  \begin{equation*}
    \xymatrix{
      {\F_2} \ar[r] & {H_*(\CPlll)} \ar[r]^-{\Omega\omega_*} &
      {H_*(Q\CP^\infty_+)} \ar[r]^-{\Omega\partial_*} & {H_*(\Omega
        QS^0)} \ar[r] & {\F_2}
    }
  \end{equation*}
  is exact.  In particular $\Omega\omega_*$ induces an isomorphism
  \begin{equation*}
    \xymatrix{ 
      {H_*(\CPlll)} \ar[r]^-{\isom} &
      {H_*(Q\CP^\infty_+)\modd\Omega\partial_*}
    }
  \end{equation*}
\end{thm}
The short exact sequence in Theorem~\ref{thm-3} determines
$H_*(\CPlll;\F_2)$ completely.  Indeed, as part of the proof of
Theorem~\ref{thm-3} we determine $H_*(\Omega QS^0)$ and the induced
map $\Omega\trf_*$ in homology.  The result can be summarised by the
diagram
\begin{equation}
  \label{eq:8}
  \xymatrix{
    {QH_*(Q\CP^\infty_+)} \ar[r]^-{Q(\Omega\trf_*)}
    \ar[d]^-{\isom} &
    {QH_*(\Omega QS^0)} \ar[d]^-{\isom} \\
    {PH_*(Q\Sigma\CP^\infty_+)} \ar[r]^-{P(\trf_*)}
    &
    {PH_*(QS^0)}
  }
\end{equation}
where the vertical isomorphisms in \eqref{eq:8} are the homology
suspensions.  The homology $H_*(\Omega_0 QS^0)$ of the basepoint
component of $\Omega QS^0$ is a \emph{divided power algebra}, i.e.\ 
its dual is a primitively generated polynomial algebra.

For odd primes $p$ our results are less precise in that
$H_*(\CPlll;\F_p)$ is only determined up to algebra isomorphism.  The
main technical theorem is the following
\begin{thm}
  \label{thm-4}
  For odd primes $p$, the homology suspension
  \begin{equation*}
    \sigma_*:  QH_*(\CPlll) \to PH_*(\CPll)
  \end{equation*}
  is surjective.
\end{thm}
Proving Theorem~\ref{thm-4} is the most difficult part of the paper.
It uses that $\sigma_*$ commutes with the homology operations
$\beta^\epsilon Q^s$.
\begin{cor}
  \label{cor-5}
  Let $Y\subseteq H_*(\CPlll)$ be a subset such that $\sigma_*(Y)$ is
  a basis of $PH_*(\CPll)$.  Then $H_*(\CPlll)$ is the free
  commutative algebra on the set
  \begin{equation*}
    Y \cup \{ \beta Q^s y \vert \text{$y\in Y^-$, $\deg(y) = 2s-1$}\}
  \end{equation*}
\end{cor}
Corollary~\ref{cor-5} is a formal consequence of Theorem~\ref{thm-4}
and the fact that the Hopf algebra $H_*(\CPll)$ is primitively
generated.  The proof uses the ``Kudo transgression theorem'', cf.\ 
\cite{CLM}, Theorem~1.1(7): If $\deg(y) = 2s-1$, then in the
Leray-Serre spectral sequence we have that $\sigma_*(y)$ transgresses
to $y$ and that $(\sigma_*y)^{p-1} \tensor y$ transgresses to $-\beta
Q^s y$.

\subsection{Acknowledgements}  This calculation is part of my
phd-project at the University of Aarhus.  It is a great pleasure to
thank my thesis advisor Ib Madsen for his help and encouragement
during my years as a graduate student.

\section{Recollections}\label{sec:recollections}

In this introductory section we collect the results we need later in
the paper.  We start by recalling some important results on the
structure of Hopf algebras from \cite{MM} and proceed to review the
functor $\Cotor$ and the closely related Eilenberg-Moore spectral
sequence, cf.\ \cite{EM}, \cite{SM}.

\subsection{Hopf algebras}
Here and elsewhere, the field $\F_p$ with $p$ elements is the ground
field, and $\tensor = \tensor_{\F_p}$.  Until further notice, $p$ is
assumed odd.  Algebras and coalgebras are as in \cite{MM} and in
particular they always have units resp.\ counits.

\begin{defn}\label{defn:1}
  When $A$ is a coalgebra and $M_A, {_A N}$ are $A$-comodules with
  structure maps $\Delta_M: M \to M \tensor A$ and $\Delta_N: N\to A
  \tensor N$, the cotensor product is defined by the exact sequence
  \begin{equation*}
    \xymatrix{
      {0} \ar[r] & {M \cotensor_A N} \ar[r] & {M\tensor N} \ar[r] &
      {M\tensor A \tensor N}
    }
  \end{equation*}
  where the right-hand morphism is ${\Delta_M \tensor N - M \tensor
    \Delta_N}$.  The functors $M\cotensor_A -$ and $-\cotensor_A N$
  are left exact functors from $A$-comodules to $\F_p$-vectorspaces in
  general, and to $A$-comodules when $A$ is cocommutative.
\end{defn}

\begin{defn}
  \label{defn:2}
  For a morphism $f: A\to B$ of Hopf algebras, define the
  \emph{kernel} and \emph{cokernel}
  \begin{equation*}
    A\modd f = A\cotensor_B k, \quad
    B\moddd f = B \tensor_A k
  \end{equation*}
\end{defn}
A priori, the kernel and cokernel are vectorspaces, but when $A$ and
$B$ are commutative and cocommutative, they become Hopf algebras and
are the kernel and cokernel in the categorical sense.  Hopf algebras
that are both commutative and cocommutative are called \emph{abelian},
and the category of those is an abelian category (this essentially
follows from \cite[Prop.\ 4.9]{MM}).

The Hopf algebras appearing in this paper will (except for $R$ and
$\RR$ defined below) be of the form $A = H_*(X;\F_p)$ for $X$ an
infinite loop space.  Such Hopf algebras will always be abelian.  We
will often have that $H_*(X,\F_p)$ is of finite type, and in this case
$H^*(X;\F_p)$ will also be a Hopf algebra.  However, if $\pi_0(X)$ is
infinite, $H^*(X;\F_p)$ will not be a Hopf algebra (e.g.\ $X=QS^0$
with $\pi_0X = \Z$).  Usually it will then be the case that the
basepoint component $X_0\subseteq X$ is of finite type, and thus we
can consider $H^*(X_0;\F_p)$.  Hopf algebras $A$ with $A_i = 0$ for
$i<0$ and $A_0 = \F_p$ are called \emph{connected}.  In general we
will have a natural splitting of Hopf algebras $H_*(X) =
H_*(X_0)\tensor \F_p\{\pi_0X\}$ where $\F_p\{\pi_0X\} = H_0(X)$ is the
group algebra.

\begin{defn}
  \label{defn:3}
  For an algebra $A$ with augmentation $\epsilon$, let $IA =
  \Ker(\epsilon: A \to k)$ and dually for a coalgebra $A$ with
  augmentation $\eta$, let $JA=\Cok(\eta:k\to A)$.  Let $Q$ and $P$ be
  the functors defined by the exact sequences
  \begin{equation*}
    \xymatrix{
      {IA\tensor IA} \ar[r]^-{\phi} & IA\ar[r]& QA \ar[r] & 0
    }
  \end{equation*}
  and
  \begin{equation*}
    \xymatrix{
      0 \ar[r] & PA \ar[r] & JA\ar[r]^-{\Delta} & {JA\tensor JA}
    }
  \end{equation*}
\end{defn}

$P$ and $Q$ satisfies $P(A\tensor B) = PA \oplus PB$ and $Q(A\tensor
B) = QA\tensor QB$, and as functors from abelian Hopf algebras to
vectorspaces, $Q$ is right exact and $P$ is left exact (\cite[Prop
4.10]{MM}).  When $A$ is connected, $PA \subseteq A$ is the subset
consisting of elements $x$ satisfying $\Delta x = x\tensor 1 +
1\tensor x$.  If $A = \F_p\{\pi_0X\}$ is a group algebra, then $PA =
0$.

The functors $P$ and $Q$ are related by the short exact sequence of
\cite[Thm.\ 4.23]{MM}:
\begin{thm}
  \label{thm:mmses}
  For an abelian Hopf algebra $A$, let $\xi: A\to A$ be the Frobenius
  map $x\mapsto x^p$ and let $\lambda: A\to A$ be the dual of $\xi:
  A^* \to A^*$.  Let $\xi A\subseteq A$ be the image of $\xi$ and let
  $A\to\lambda A$ be the coimage of $\lambda$.  Then there is the
  following natural exact sequence
  \begin{equation}
%    \label{eq:mmses}
    \xymatrix{
      0\ar[r]&P\xi A\ar[r]& PA\ar[r]&QA\ar[r]& {Q\lambda A}\ar[r]& 0
    }
  \end{equation}
  In particular $PA\to QA$ is an isomorphism except possibly in
  degrees $\equiv 0\pmod{2p}$ if $p>2$.  For $p=2$ it is an
  isomorphism in odd degrees.
\end{thm}

Finally, we recall Borel's structure theorem (\cite[Theorem 7.11]{MM})
\begin{thm}
  \label{thm:borel-structure}
  Any connected abelian Hopf algebra $A$ is isomorphic as an algebra
  to a tensor product of algebras of the form $E[x]$, $\F_p[x]$ and
  $\F_p[x]/(x^{p^n})$, $n\geq 1$.
\end{thm}

\begin{cor}\label{cor:fricoalg}
  A connected abelian Hopf algebra $A$ is isomorphic as an algebra to
  a polynomial algebra if and only if $\xi: A\to A$ is injective.
  Dually if $A$ is of finite type, $A^*$ is polynomial if and only if
  $\lambda: A\to A$ is surjective.\qed
\end{cor}

\subsection{The functor $\Cotor$}  When $A$ is a coalgebra and $B$ and
$C$ are left resp.\ right $A$-comodules, the functor
\begin{displaymath}
  \Cotor^A(B,C)
\end{displaymath}
is defined as the right derived functor of the cotensor product
$\cotensor_A$.  To be explicit (and to fix grading conventions),
choose an injective resolution $0 \to B \to I_0 \to I_{-1} \to \dots$
of $B$ in the category of right $A$-comodules and set
\begin{displaymath}
  \Cotor^A_n(B,C) = H_n(I_* \cotensor_A C)
\end{displaymath}
When $A$, $B$ and $C$ are in the graded category, $\Cotor$ gets an
inner grading and is thus bigraded with $\Cotor^A_{n,m}(B,C) =
(\Cotor^A_n(B,C))_m$.  When $A,B,C$ are all positively graded,
$\Cotor$ is concentrated in the second quadrant.

When $A$, $B$ and $C$ are of finite type over a field, this functor is
dual to the more common $\Tor$:
\begin{displaymath}
  \Cotor^A(B,C) = \left(\Tor^{A^*}(B^*,C^*)\right)^*
\end{displaymath}
This follows immediately from the duality between $\cotensor_A$ and
$\tensor_{A^*}$.

We shall consider $\Cotor$ as a functor from diagrams of cocommutative
coalgebras
\begin{equation*}
  \mathscr{S} = \left\{
    \begin{array}{c}
      \xymatrix{
        & B\ar[d]\\
        C\ar[r] & A\\
      }
    \end{array}
  \right\}
\end{equation*}
to coalgebras.  The external product is an isomorphism (see
\cite[Theorem~3.1, p.\ 209]{CE})
\begin{displaymath}
  \Cotor^A(B,C)\tensor\Cotor^{A'}(B',C') \to
  \Cotor^{A\tensor A'}(B\tensor B',C\tensor C')
\end{displaymath}
and under this isomorphism the comultiplication in $\Cotor^A(B,C)$ is
given by the comultiplication $\Delta: \mathscr{S}\to
\mathscr{S}\tensor \mathscr{S}$ in the diagram $\mathscr{S}$.

Dually, when $\mathscr{S}$ is a diagram of Hopf algebras,
$\Cotor^A(B,C)$ is a Hopf algebra with multiplication induced by the
multiplication $\phi: \mathscr{S}\tensor\mathscr{S} \to \mathscr{S}$
of the diagram $\mathscr{S}$.

Later we will need the structure of $\Cotor^A(B,\F_p)$ where $\F_p$
denotes the trivial Hopf algebra and $f: B\to A$ is a morphism of Hopf
algebras.  From the change of rings spectral sequence and
\cite[Theorem 4.9]{MM} we get
\begin{prop}\label{prop:Cotorsplit}
  For a map $f: B\to A$ of Hopf algebras, there is a natural
  isomorphism of Hopf algebras
  \begin{displaymath}
    \xymatrix{
      {\Cotor^A(B,\F_p)} \ar[r]^-{\isom} & {B\modd f} \tensor
      {\Cotor^{A\moddd f}(\F_p,\F_p)}
    }
  \end{displaymath}\qed
\end{prop}

To complete the description of $\Cotor^A(B,\F_p)$ we need to calculate
the Hopf algebra $\Cotor^A(\F_p,\F_p)$.  This is easily done by
applying Borel's structure theorem to the dual algebra $A^*$ and using
Lemma~\ref{lem:Torstruktur} below.  The Hopf algebra $\Gamma[x]$ is
dual to a polynomial algebra: $\Gamma[x] = (k[x^*])^*$ and $s^{-\nu}$
denotes bigraded desuspension: $(s^{-\nu}V)_{-\nu,n} = V_n$ for a
singly graded object $V$.
\begin{lem}\label{lem:Torstruktur}
  The following isomorphisms hold as Hopf algebras
  \begin{align*}
    \Tor^{E[x]}(\F_p,\F_p) & = \Gamma[s^{-1}x]\\
    \Tor^{\F_p[x]}(\F_p,\F_p) & = E[s^{-1}x]\\
    \Tor^{\F_p[x]/(x^{p^n})}(\F_p,\F_p) & = E[s^{-1}x]\tensor \Gamma[s^{-2}
    x^{p^n}]
  \end{align*}\qed
\end{lem}
By the duality between $\Tor$ and $\Cotor$ we obtain the Hopf algebra
structure of $\Cotor^A(\F_p,\F_p)$ in terms of a set of generators of the
dual algebra $A^*$.
\begin{cor}\label{lem:Cotorstruktur}
  For any connected Hopf algebra $A$ of finite type,
  $\Cotor^A(\F_p,\F_p)$ is freely generated by the primitive elements
  in $\Cotor^A_{-1,*}(\F_p,\F_p)$ and $\Cotor^A_{-2,*}$.  Choosing
  generators of $A^*$ (according to Borel's structure theorem), the
  generators of $\Cotor^A(\F_p,\F_p)$ are in bidegrees
  \begin{align*}
    (-1,k) & \quad\text{for $x\in A_k^*$ an odd generator}\\
    (-1,k) & \quad\text{for $x\in A_k^*$ an even generator}\\
    (-2,p^m k) & \quad\text{for $x\in A_k^*$ an even generator of
    height $p^m$}
  \end{align*}
  The primitive elements of $\Cotor^A(k,k)$ are in bidegrees
  \begin{align*}
    p^n(-1,k) & \quad\text{for $x\in A_k^*$ an odd generator}\\
    (-1,k) & \quad\text{for $x\in A_k^*$ an even generator}\\
    p^n(-2,p^m k) & \quad\text{for $x\in A_k^*$ an even generator of
    height $p^m$}
  \end{align*}
  \qed
\end{cor}
More functorially, one defines for $p>2$ the functor $\Hat PA =
P\Cotor^A_{-2,*}(\F_p,\F_p)$.  Then the result in
Corollary~\ref{lem:Cotorstruktur} is that $\Cotor^A(\F_p,\F_p) \isom
S[s^{-1}PA]\tensor S[s^{-2}\Hat PA]$ combined with the facts that
$\Hat Q(\F_p[x]) = \Hat Q(E[x]) = 0$ and that $\Hat
Q(\F_p[x]/(x^{p^n})) = \F_p.\{x^{p^n}\}$, where $\Hat QA = (\Hat
PA^*)^*$.

In particular, the only primitive elements of \emph{odd} total degree
are in bidegrees $(-1,k)$ for even generators $x\in A_k^*$.

Finally, we shall need a criterion for left exactness of the functor
$Q$, namely
\begin{prop}
  \label{prop:Q-inj}
  Let
  \begin{equation*}
    k \to A\to B\to C\to k
  \end{equation*}
  be a short exact sequence of abelian Hopf algebras.  If $C$ is a
  free commutative algebra, then the sequence
  \begin{equation*}
    0 \to QA\to QB\to QC\to 0
  \end{equation*}
  is short exact.
\end{prop}
\begin{proof}
  Since $C$ is free, we may split $B\to C$ with a map of algebras.
  Thus $B \isom A\tensor C$ as an algebra, and $Q(B)$ depends only on
  the algebra structure of $B$.
\end{proof}
A peculiar consequence of Corollary~\ref{cor:fricoalg} is that if $A$
is a Hopf algebra that is free as an algebra, then any Hopf subalgebra
of $A$ is also free as an algebra.

\subsection{The spectral sequence}  In this section, we recall the
spectral sequence of \cite{EM} and some of its properties.

We consider \emph{homotopy cartesian squares}
\begin{displaymath}
  \mathscr{C}= 
  \left\{
    \begin{array}{c}
      \xymatrix{
        F\ar[r]\ar[d] & E\ar[d]\\
        X\ar[r] & B\\
      }
    \end{array}
  \right\}
\end{displaymath}
of connected spaces, and with $B$ simply connected (homotopy cartesian
means that $F \simeq \holim(X\to B \leftarrow E)$.  One can always find a
model that is a \emph{fibre square}, i.e. where $E\to B$ is a
fibration, and $F\to X$ is the pullback fibration).  In the following,
$H_*$ denotes $H_*(-;\F_p)$.
\begin{defn}
  The Eilenberg-Moore spectral sequence $E^r$ is a functor from fibre
  squares $\mathscr{C}$ as above to spectral sequences of coalgebras.
  It has
  \begin{displaymath}
    E^2 = \Cotor^{H_*(B)}(H_*(E),H_*(X))
  \end{displaymath}
  and converges as coalgebra to $H_*F$.
\end{defn}

\begin{thm}[{\cite[Proposition 16.4]{EM}}]\label{thm:externalp}
  The external product induces an isomorphism
  \begin{displaymath}
    E^r(\mathscr{C}) \tensor E^r(\mathscr{C}') \to
    E^r(\mathscr{C}\times \mathscr{C}')
  \end{displaymath}
  Under this isomorphism, the coalgebra structure is induced by the
  diagonal $\Delta: \CC\to \CC\times\CC$.\qed
\end{thm}

Dually, when $\CC$ is a diagram of $H$-spaces and $H$-maps (here
meaning maps commuting \emph{strictly} with the multiplication such as
loop spaces and loop maps), there is a multiplication $m: \CC\times\CC
\to \CC$ inducing a multiplication $\phi = m_*: E^r(\CC)\tensor
E^r(\CC) \to E^r(\CC)$.  In this case, the spectral sequence is one of
Hopf algebras.  Furthermore it is clear that on the $E^2$-term, the
Hopf algebra structure is the same as the one on $\Cotor$ described
above.

\subsection{The loop suspension}  We shall use the spectral
sequence only in the case when $X$ is a point.  This corresponds to a
fibration sequence
\begin{displaymath}
  F\to E\to B
\end{displaymath}
and the spectral sequence computes homology of
the fibre.  When $E$ is also a point, we have the path-loop fibration
sequence
\begin{displaymath}
  \Omega X \to *\to X
\end{displaymath}
In this case, the line
\begin{displaymath}
  E^2_{0,*} = \Cotor^{H_*(X)}_{0,*}(\F_p,\F_p) = \F_p\cotensor_{H_*(X)} \F_p = \F_p
\end{displaymath}
is concentrated in degree 0 and hence there is a ``secondary edge
homomorphism''
\begin{equation}\label{eq:susp}
  H_*(\Omega X) \to E^\infty_{-1,*} \hookrightarrow E^2_{-1,*} \isom
  PH_*X
\end{equation}
\begin{prop}[{\cite[Proposition 4.5]{S}}]\label{prop:susp}
  The morphism in \eqref{eq:susp} is the loop suspension
  \begin{displaymath}
    \sigma_*: \Ind H_*(\Omega X) \to \Prim H_*X
  \end{displaymath}
\end{prop}\qed

We shall also need
\begin{lem}\label{lem:diffHopf}
  Let $C_*$ be a connected differential graded Hopf algebra. If $x$ is
  an element of minimal degree with $dx \neq 0$, then $x$ is
  indecomposable and $dx$ is primitive.
\end{lem}
\begin{proof}
  Immediate from the Leibniz rules for product and coproduct.
\end{proof}
\begin{cor}\label{cor:mindiff}
  Minimal differentials in the spectral sequence of a path-loop
  fibration correspond to minimal elements in the cokernel of
  $\sigma_*$.
\end{cor}
\begin{proof}
  Since $dx$ is primitive and not in $E^2_{-1,*}$ it is of even total
  degree and $x$ is of odd total degree.  By
  Corollary~\ref{lem:Cotorstruktur}, the only odd dimensional indecomposable
  elements are in $E^2_{-1,*}$ and the result follows.
\end{proof}

\section{Unstable $R$-modules}
\label{sec:unstable-r-modules-1}

As sketched in the introduction, homology of an infinite loop space
has homology operations $\beta^\epsilon Q^s$.  In this section we
recall the precise definitions and explain how to express homology of
$QX$ as a free algebra on certain iterated operations on the homology
of $X$.  We follow the notation from \cite{CLM}.  In this section we
consider only $p>2$.  Small changes, which we recall later, are needed
for $p=2$.

We define several categories of graded vectorspaces with a set of
linear transformations $\{\beta^\epsilon Q^s \mid \epsilon\in\{0,1\},
s\in\Z_{\geq\epsilon}\}$ of degree $2s(p-1)-\epsilon$.

\begin{equation}\label{eq:1}
  \xymatrix{
    {\text{$Q$-unstable $R$-modules}} \ar[r]\ar[d] &
    {\text{$Q$-unstable $\RR$-modules}} \ar[d]\\
    {\text{unstable $R$-modules}} \ar[r] &
    {\text{unstable $\RR$-modules}}\ar[d] \\
    & {\text{$\RR$-modules}}\ar[d]\\
    & {\text{graded vectorspaces}}\\
  } 
\end{equation}
Here, $\RR$ is the free non-commutative algebra on the set
$\{\beta^\epsilon Q^s \mid \epsilon\in\{0,1\},
s\in\Z_{\geq\epsilon}\}$, and the various entries in \eqref{eq:1}
differ in what relations the action of the operations $\beta^\epsilon
Q^s$ are assumed to satisfy.  It is the left part of the diagram that
is geometrically relevant, since the homology of an infinite loop
space $X$ is naturally an unstable $R$-modules, and so is the
vectorspace of primitive elements $PH_*(X)$.  The space of
indecomposable elements $QH_*(X)$ is naturally a $Q$-unstable
$R$-module.

All of the above forgetful functors to graded vectorspaces have left
adjoint ``free'' functors.  From $\RR$-modules it is the functor
$V\mapsto \RR\tensor V$, and the other four are quotients thereof.

In \ref{sec:algebras-rr-r-1}, we define the algebras $\RR$ and $R$ and
the four categories of unstable modules.  In \ref{sec:free-functors}
we construct the four adjoint functors $\DD$, $\DD'$, $D$ and $D'$.
Finally, in \ref{sec:homology-qx} we recall the computation of
$H_*(QX)$ in terms of $H_*(X)$.  It should be noted that the algebra
$\RR$ and the related categories are needed only in the proof of
Theorem~\ref{thm:4.3}.  It is $R$ that is geometrically relevant but
$\RR$ has the property that a submodule of a free $\RR$-module is
again free and similarly for submodules of free ($Q$-)unstable
$\RR$-modules.  This makes $\RR$ simpler from the viewpoint of
homological algebra.

\subsection{Araki-Kudo-Dyer-Lashof operations}
\label{sec:kudo-araki-dyer}

Recall that an \emph{infinite loop space} is a sequence $E_0, E_1,
\dots$ of spaces and homotopy equivalences $\Omega E_{i+1} \to E_i$.
One thinks of $E_0$ as the ``underlying space'' of the infinite loop
space.  In particular, $E_0 = \Omega^2E_2$ is a homotopy commutative
$H$-space.  Thus, as mentioned in the introduction, $H_*(E_0)$ is a
commutative algebra under the Pontrjagin product.  Furthermore
$H_*(E_0)$ naturally carries a set of linear transformations
$\beta^\epsilon Q^s$, $\epsilon\in\{0,1\}, s\in\Z_{\geq{\epsilon}}$.
These linear transformations are commonly called \emph{Dyer-Lashof
  operations} (or \emph{Araki-Kudo} operations) and are operations
\begin{displaymath}
  \beta^\epsilon Q^s: H_n(E_0) \to H_{n+2s(p-1)-\epsilon}(E_0)
\end{displaymath}
natural with respect to infinite loop maps. They measure the failure
of chain level commutativity of the Pontrjagin product.

They satisfy a number of relations that makes $H_*(E_0)$ an
\emph{unstable $R$-module}, the notion of which is defined below.

\subsection{The algebras $\RR$ and $R$ and categories of unstable
  modules}
\label{sec:algebras-rr-r-1}

\begin{defn}
  Let $\RR$ be the free (non-commutative) algebra generated by symbols
  \begin{displaymath}
    \beta^\epsilon Q^s,\quad \epsilon\in\{0,1\},
    s\in\Z_{\geq\epsilon}.
  \end{displaymath}
  and write $\beta Q^s = \beta^1 Q^s$ and $Q^s = \beta^0 Q^s$.  $\RR$
  is a graded algebra with
  \begin{displaymath}
    \deg(\beta^\epsilon Q^s) = 2s(p-1)-\epsilon
  \end{displaymath}
  It will occasionally be convenient to consider $\RR$ as a
  \emph{bigraded} algebra with gradings
  \begin{equation*}
    \deg_Q(\beta^\epsilon Q^s) = 2s(p-1),\quad
    \deg_\beta(\beta^\epsilon Q^s) = -\epsilon
  \end{equation*}
  $\RR$ is a cocommutative Hopf algebra with comultiplication
  \begin{displaymath}
    \Delta (\beta^\epsilon Q^s) = \sum_{\substack{\epsilon_1 + \epsilon_2 =
    \epsilon\\ s_1 + s_2 = s}} \beta^{\epsilon_1} Q^{s_1}\tensor
    \beta^{\epsilon_2} Q^{s_2}    
  \end{displaymath}
\end{defn}
\begin{rem}
  $\RR$ is a Hopf algebra in the sense of~\cite{MM}, i.e.\ a monoid
  object in the category of cocommutative coalgebras.  Notice however
  that $\RR$ is not a group object, since $Q^0$ is not invertible.
\end{rem}
\begin{defn}
  An $\RR$-module is called \emph{unstable}, if
  \begin{equation}\label{eq:unstability}
    \beta^\epsilon Q^s x = 0\quad\text{whenever $2s-\epsilon <
    \deg(x)$}
  \end{equation}
  It is called \emph{$Q$-unstable} if furthermore
  \begin{equation}
    \label{eq:2}
    Q^s x = 0\quad\text{whenever $2s=\deg(x)$}
  \end{equation}
\end{defn}
On homology of an infinite loop space we also have the relation
\begin{equation}
  \label{eq:6}
    Q^s x = x^p\quad\text{whenever $2s = \deg(x)$}
\end{equation}

For an infinite loop space $X$, $H_*(X)$ is naturally an unstable
$\RR$-module and $QH_*(X)$ is $Q$-unstable because of~\eqref{eq:6}.
However, the ideal in $\RR$ of elements with universally trivial
action is nonzero, and hence the action of $\RR$ on $H_*X$ factors
through a quotient of $\RR$.  This quotient is the \emph{Dyer-Lashof
  algebra} $R$.

\begin{defn}\label{defn:Adem}
  For each $r,s\in\N$ and $\epsilon \in\{0,1\}$ with $r > ps$, define
  elements in $\RR$
  \begin{displaymath}
    \AA^{(\epsilon,r,0,s)} = \beta^\epsilon Q^rQ^s - \left(
    \sum_{i=0}^{r+s} (-1)^{r+i} (pi-r, r-(p-1)s -i -1) \beta^\epsilon
    Q^{r+s-i} Q^i\right)
  \end{displaymath}
  For $r\geq ps$ define elements
  \begin{displaymath}
    \begin{split}
      \AA^{(0,r,1,s)} = Q^r\beta Q^s - \bigg(
        \sum_{i=0}^{r+s} (-1)^{r+i} (pi-r, r-(p-1)s -i) \beta
        Q^{r+s-i} Q^i \\- \sum_{i=0}^{r+s} (-1)^{r+i} (pi-r-1, r-(p-1)s -i)
        Q^{r+s-i} \beta Q^i \bigg)
    \end{split}
  \end{displaymath}
  and
  \begin{displaymath}
    \AA^{(1,r,1,s)} = \beta Q^r\beta Q^s - \bigg(- \sum_{i=0}^{r+s}
    (-1)^{r+i} (pi-r-1, r-(p-1)s -i) \beta Q^{r+s-i} \beta Q^i \bigg)
  \end{displaymath}
  where $(i,j) = (i+j)!/(i!j!)$. These elements are the socalled
  \emph{Adem relations}.
\end{defn}

Let $\AA\subseteq \RR$ be the $\F_p$-span of all Adem elements.  This
is a bigraded subspace of $\RR$.  Let $\langle\AA\rangle\subseteq \RR$
be the two-sided ideal generated by $\AA$.  Let $\JJ\subseteq \RR$ be
the two-sided ideal (or equivalently the left ideal) generated by the
relations \eqref{eq:unstability} (for $x\in\RR$). $\JJ$ is the
smallest ideal such that $\RR/\JJ$ is unstable as a left $\RR$-module.
\begin{defn}
  The \emph{Dyer-Lashof algebra} is the quotient
  \begin{displaymath}
    R = \RR/(\langle \AA\rangle + \JJ)
  \end{displaymath}
\end{defn}

The action of $\AA$ and hence $\langle\AA\rangle$ on homology of
infinite loop spaces is trivial by results from \cite{CLM}, dual to
Adem's result for the Steenrod algebra.  So is the action of $\JJ$.
Hence $H_*(X)$ is an $R$-module when $X$ is an infinite loop space.
Conversely (\cite{AK}, \cite{DL}) the map $R\to H_*(QS^0)$ induced by
acting on the zero-dimensional class $\iota$, corresponding to the
non-basepoint of $S^0$, is an injection, so there are no further
relations.

The set of all products of generators form a vector space basis of
$\RR$.  To have an explicit basis for $R$, we recall the notion of
admissible monomials, \cite[p.\ 16]{CLM}.

A sequence
\begin{equation*}
  I = (\epsilon_1, s_1, \dots, \epsilon_k, s_k)
\end{equation*}
of integers $\epsilon_i\in\{0,1\}$ and $s_i\in\Z_{\geq \epsilon_i}$
determines the iterated homology operation
\begin{equation*}
  Q^I = \beta^{\epsilon_1}Q^{s_1}
  \dots\beta^{\epsilon_k}Q^{s_k}\in\RR
\end{equation*}
This sequence is called \emph{admissible} if for all $i = 2,\dots,
k$,
\begin{equation}\label{eq:30}
  s_i \leq ps_{i-1} - \epsilon_{i-1}
\end{equation}
The corresponding iterated homology operations $Q^I\in\RR$ are called
\emph{admissible monomials}.  The \emph{length} and \emph{excess} of
$I$ are
\begin{equation*}
  \ell(I) = k,\quad e(I) = 2s_1 - \epsilon_1 - \sum_{j=2}^k [2s_j(p-1) - \epsilon_j]
\end{equation*}
Furthermore, define
\begin{equation*}
  b(I) = \epsilon_1
\end{equation*}

Using the Adem relations one may rewrite an arbitrary element of $R$
as a linear combination of admissible monomials.  Applying Adem
relations does not raise the excess.

There is a natural quotient map $\RR\to R$.  Thus $R$-modules are also
$\RR$-modules.
\begin{defn}
  An $R$-module is called \emph{unstable}, respectively
  \emph{$Q$-unstable}, if it is so as an $\RR$-module.
\end{defn}

\subsection{Free functors}
\label{sec:free-functors}

\begin{defn}
  For a graded vectorspace $V$ we define $\DD V$ to be the quotient of
  $\RR\tensor V$ by the relations \eqref{eq:unstability} and $\DD'V$
  to be the quotient of $\DD V$ by the relations \eqref{eq:2}.  Define
  also
  \begin{equation*}
    DV = R\tensor_\RR \DD V,\quad
    D'V = R\tensor_\RR \DD' V
  \end{equation*}
\end{defn}
The functor $\DD$ is left adjoint to the forgetful functor from
unstable $\RR$-modules to vectorspaces.  Thus $\DD V$ is the ``free
unstable $\RR$-module'' generated by $V$.  Similarly, $D$ is left
adjoint to the forgetful functor from unstable $R$-modules to graded
vectorspaces.  Analogous remarks apply to $\DD'$ and $D'$. The
functors appear in the following exact sequences, natural in $V$
\begin{align}
  \label{eq:exactseq}
  \langle \AA\rangle\tensor_\RR \DD V \to \DD V \to DV &\to 0\\
  \label{eq:exactseq2}
  \langle \AA\rangle\tensor_\RR \DD' V \to \DD' V \to D'V &\to 0
\end{align}
When $V = \F_p\iota$ for a homogeneous element $\iota$, $DV$ has basis
\begin{displaymath}
  \{Q^I \iota \mid \text{$I$ admissible, $e(I)\geq \deg(\iota)$}\}
\end{displaymath}
Together with additivity of $D$, this describes $DV$ as a
$\F_p$-vectorspace.  Since $R\isom D\F_p$ as a left $R$-module, we
also have a basis of $R$ over $\F_p$.

\subsection{Homology of $QX$}
\label{sec:homology-qx}

Here we recall the computation of $H_*(QX)$.  It can be expressed as a
functor of $H_*(X)$ which is left adjoint to a suitable forgetful
functor, forgetting the Pontrjagin product and the $R$-action,
see~\cite{CLM}.  We shall give a non-functorial description in terms
of a basis of $JH_*(X)$.

\begin{thm}
  Let $B\subseteq JH_*(X)$ be a basis consisting of homogeneous
  elements.  Then $H_*(QX)$ is the free commutative algebra on the set
  \begin{equation*}
    \mathbf{T} = \{Q^Ix\vert \text{$x\in B$, $I$ admissible, $e(I) +
    b(I) > \deg(x)$}\}
  \end{equation*}\qed
\end{thm}
\begin{cor}
  \label{cor:3.16}
  The natural map
  \begin{equation*}
    \phi_Q: D'JH_*(X) \to QH_*(QX)
  \end{equation*}
  is an isomorphism of $Q$-unstable $R$-modules.
  
  If $X$ is connected and $H_*(X)$ has trivial comultiplication (e.g.\ 
  if $X$ is a suspension), then the natural map
  \begin{equation*}
    \phi_P: DJH_*(X) \to PH_*(QX)
  \end{equation*}
  is an isomorphism of unstable $R$-modules.\qed
\end{cor}
\begin{rem}
  $QX$ is connected if and only if $X$ is connected.  More generally
  the group of components of $QX$ is determined by the short exact
  sequence
  \begin{equation*}
    0 \to \Z \to \Z[\pi_0X] \to \pi_0(QX) \to 0
  \end{equation*}
  where the first arrow is induced by the inclusion of the basepoint in
  $X$.  When $X$ is nonconnected we are sometimes only interested in
  homology of the component $Q_0X$ of the basepoint in $QX$.  This can
  be described as follows.  Let $\tau: QX \to Q_0X$ be the map that on
  a component $Q_iX$, $i\in \pi_0(QX)$ multiplies by an element of
  $Q_{-i}X$.  This defines a welldefined homotopy class of maps
  $QX\to Q_0X$ which is left inverse to the inclusion.  Then we have
  that $H_*(Q_0X)$ is the free commutative algebra on the set
  \begin{equation*}
    \mathbf{\Tilde T} = \{\tau_* Q^Ix\vert \text{$x\in B$, $I$
    admissible, $e(I) + b(I) > \deg(x)$, $\deg(Q^Ix)> 0$}\}
  \end{equation*}
\end{rem}

\section{Homological algebra of unstable modules}
\label{sec:homol-algebra-unst}

The map
\begin{equation*}
  Q(\trf_*): QH_*(Q\Sigma\CP^\infty_+) \to QH_*(Q_0S^0)
\end{equation*}
was computed in \cite[Theorem 4.5]{MMM}.  The left hand side is
$D'JH_*(\Sigma\CP^\infty_+)$ and the right hand side is $D'\F_p$.  The
starting point of our theorems is
\begin{thm}[\cite{MMM}]\label{thm:MMM}
  Let $a_s\in H_s(\Sigma\CP^\infty_+)$ be the generator, $s$ odd.  Let
  $\iota \in JH_0(S^0)$ be the generator. Then
  \begin{equation*}
    Q(\trf_*)(a_s) = 
    \begin{cases}
      (-1)^r \beta Q^r\iota & \text{$s=2r(p-1)-1$}\\
      0 & otherwise\\
    \end{cases}
  \end{equation*}
\end{thm}
\begin{proof}
  The map $\trf: \Sigma\CP^\infty_+ \to QS^0$ coincides with the
  universal $S^1$-transfer denoted $t_0$ in \cite{MMM}.  The formula
  for $Q(\trf_*)(a_s)$ in the theorem now follows from ignoring all
  decomposable terms in \cite[Theorem 4.5]{MMM}.
\end{proof}

\subsection{Main technical theorems}
\label{sec:main-techn-theor}

To state the theorems, recall from
subsection~\ref{sec:algebras-rr-r-1} that $\RR$ may be bigraded by
$\deg = \deg_Q + \deg_\beta$.  Since the Adem relations are
homogeneous with respect to $\deg_Q$ and $\deg_\beta$, there is an
induced bigrading of $R$.  If $V$ is bigraded, $\RR\tensor V$ is a
bigraded left $\RR$-module.  Since the unstability relations
\eqref{eq:unstability} are homogeneous, there is an induced
bigrading of $\DD V$.  Similarly for $\DD'V$, $DV$ and $D'V$.  Thus by
Corollary~\ref{cor:3.16} a bigrading of $JH_*(X)$ will induce a
bigrading of $QH_*(QX)$ and, for $X$ a suspension, a bigrading of
$PH_*(QX)$.  However, $H_*(QX)$ will only have $\deg_\beta$
welldefined up to multiplication with $p$ because of the unstability
relation~\eqref{eq:6}.

For bigraded modules $V$ with $\deg = \deg_Q + \deg_\beta$ as above,
we shall write $V^{i,j} = \{x\in V\mid \deg_Q(x) = i, \deg_\beta(x) =
j\}$ and $V^n = \oplus_{i+j=n} V^{i,j}$ and $V^{(n)} = \oplus_i
V^{i,n}$.  We will only consider gradings in the fourth quadrant,
i.e.\ $V^{i,j} = 0$ unless $i\geq 0$ and $j\leq 0$.  Write $V^{(-)} =
\oplus_{n<0} V^{(n)}$.

\begin{thm}\label{thm:imF}
  Bigrade $JH_*(S^0)$ by setting $\deg_\beta(\iota) = 0$ and give
  $QH_*QS^0$ the induced bigrading.  Then we have
  \begin{equation*}
    \IM(Q(\trf_*)) = QH_*(QS^0)^{(-)}
  \end{equation*}
\end{thm}
\begin{proof}
  The inclusion $\IM(Q\trf_*)\subseteq QH_*(QS^0)^{(-)}$ is
  immediate from Theorem~\ref{sec:main-techn-theor}.  The other
  inclusion follows from Lemma~\ref{lem:4.4} below.  Indeed, the
  two-sided ideal in $R$ generated by the set $\{\beta Q^s\mid s\geq
  1\}$ is spanned by operations $Q^I$ with at least one $\beta$.  By
  Lemma~\ref{lem:4.4} below, any such operation is also in the left
  ideal with the same generators, i.e.\ is a linear combination of
  elements of the form $Q^J\beta Q^s$.  In particular, any element in
  $QH_*(QS^0)^{(-)}$ is also in $\IM(Q\trf_*)$ because $Q\trf_*$ is $R$-linear.
\end{proof}
\begin{lem}\label{lem:4.4}
  The left ideal in $R$ generated by the set $\{\beta Q^s \mid
  s\geq1\}$ is also a right ideal.
\end{lem}
\begin{proof}
  Write $I\subseteq R$ for the left ideal generated by $\{\beta Q^s
  \mid s\geq1\}$.

  For $r \leq ps$, consider the Adem relation
  $\AA^{(0,ps,1,r-(p-1)s)}$:
  \begin{displaymath}
    \begin{split}
      Q^{ps}\beta Q^{r-(p-1)s} &= \beta Q^r Q^s \\
      + &\sum_{i>s} \lambda_i
      \beta Q^{r+s-i}Q^i\\
      + &\text{ terms of form $Q^{r+s-i}\beta Q^i$}
    \end{split}
  \end{displaymath}
  where we have singled out the term in the Adem relation
  corresponding to $i=s$, and where the $\lambda_i\in k$ are certain
  binomial coefficients. This shows that in the left $R$-module $R/I$
  we can write $\beta Q^r Q^s$ as a linear combination of $\beta Q^a
  Q^b$ with $a < r$. In particular, $\beta Q^1 Q^s = 0 \in R/I$ and
  by induction $\beta Q^r Q^s = 0\in R/I$.
  
  Thus we have $\beta Q^r Q^s \in I$ whenever $\beta Q^rQ^s$ is
  admissible. Since a nonadmissible $\beta Q^r Q^s$ is a linear
  combination of admissible ones, we have $\beta Q^r Q^s\in I$ for
  any $r,s$. This shows that $I$ is invariant under right
  multiplication with $Q^s$. Since it is obviously invariant under
  right multiplication with $\beta Q^s$ it follows that $I$ is a
  right ideal.
\end{proof}

The kernel of $Q\trf_*$ is harder to determine explicitly.  The
partial information contained in Theorem~\ref{thm:4.3} below suffices
for the calculation.

Notice that for any $\RR$-module $V$, the augmentation of $\RR$ gives
a natural quotient map $V \to \F_p\tensor_\RR V$ identifying
$\F_p\tensor_\RR V$ with the quotient of $V$ by the relations
$\beta^\epsilon Q^s x = 0$ for $x\in V, \epsilon\in\{0,1\}, s \geq
\epsilon$.  The functor $\F_p\tensor_\RR -$ agrees with the functor
$\F_p\tensor_R -$ on $R$-modules.  Thus the vectorspace $\F_p
\tensor_\RR V = \F_p\tensor_R V$ measures the dimensions of a minimal
set of $R$-module generators of an unstable $R$-module $V$.

In the next theorem, $a_s \in JH_s(\Sigma\CP^\infty_+)$ denotes the
generator for $s$ odd.
\begin{thm}
  \label{thm:4.3}
  Bigrade $JH_*(\Sigma\CP^\infty_+)$ by concentrating it in
  $\deg_\beta = -1$ and give $QH_*(Q\Sigma\CP^\infty_+)$ the induced
  bigrading.  Then the bigraded vectorspace
  \begin{equation*}
    \F_p\tensor_R \Ker(Q\trf_*) = \F_p\tensor_R
    Q(H_*(Q\Sigma\CP^\infty_+)\modd\trf_*)
  \end{equation*}
  is concentrated in bidegrees $\deg_\beta = -1$ and $\deg_\beta =
  -2$.  In particular $\Ker(Q\trf_*)$ is generated as an $R$-module by
  the elements $a_s\in\Ker(Q\trf_*)$ with $s\not\equiv
  -1\pmod{2(p-1)}$ together with elements of degree $\equiv -1$ and
  $\equiv -2 \pmod{2(p-1)}$.
\end{thm}
\begin{proof}
  The equality $\Ker(Q\trf_*) =
  Q(H_*(Q\Sigma\CP^\infty_+)\modd\trf_*)$ in the theorem follows from
  Proposition~\ref{prop:Q-inj} because $H_*(QS^0)$ is a free
  commutative algebra.
  
  The last statement of the theorem follows from the first.  Indeed
  the elements $Q^I a_s$ are all in the kernel of $Q(\trf_*)$ when
  $s\not\equiv -1 \pmod{2(p-1)}$ because $a_s$ is in the kernel.
  These elements give rise to one ``tautological'' element $a_s \in
  F_p \tensor_R \Ker(Q\trf_*)$.  On the span of the $Q^I a_s$ with $s
  \equiv -1 \pmod{2(p-1)}$ the claim about degrees of generators
  follows since on these elements $\deg \equiv \deg_\beta
  \pmod{2(p-1)}$.  Thus we need only prove the first statement of the
  theorem.

  We have the short exact sequence of $Q$-unstable $R$-modules
  \begin{equation}
    \label{eq:19}
    \xymatrix{
      0\ar[r]& {\Ker(Q\trf_*)}\ar[r] &
      {QH_*(Q\Sigma\CP^\infty_+)}\ar[r]^-{Q\trf_*} &
      {QH_*(QS^0)^{(-)}} \ar[r]&0
    }
  \end{equation}
  If were to apply the functor $\F_p\tensor_R -$ from $R$-modules to
  vectorspaces, we would get a long exact sequence involving
  $\Tor^R_*(\F_p,-)$, and a determination of the map induced by
  $Q\trf_*$ in $\Tor_1$ would give the result.  This is more or less
  what we do, except that it is technically more convenient to replace
  the functor $\F_p\tensor_R -$ by $\F_p\tensor_\RR -$ and to replace
  $\Tor$ by a suitable functor taking unstability into account.  We
  proceed to make these ideas precise.

  The category of $Q$-unstable $\RR$-modules is abelian and has enough
  projectives.  The functor $\F_p\tensor_\RR-$ from $Q$-unstable
  $\RR$-modules is right exact, hence the left derived functors
  $L_r(\F_p\tensor_\RR-)$ are defined.  These are unstable versions of
  $\Tor^\RR_r(\F_p,-)$.  For brevity, let us write $T^\RR_1(\F_p,-) =
  L_1(\F_p\tensor_\RR -)$.

  With these definitions, applying the functor $\F_p\tensor_\RR-$ to the
  sequence \eqref{eq:19} induces the exact sequence
  \begin{equation}
    \label{eq:20}
    \xymatrix{
      0\ar[r]& {\Cok(T^\RR_1(\F_p,Q\trf_*))} \ar[r] &
      {\F_p\tensor_R\Ker(Q\trf_*)} \ar[r]& {\Ker(\F_p\tensor_RQ\trf_*)}
      \ar[r]& 0
    }
  \end{equation}
  \emph{Claim 1.} The elements $a_s\in\Ker(Q\trf_*)$ with $s\not
  \equiv -1 \pmod{2(p-1)}$ maps in \eqref{eq:20} to a generating set
  in $\Ker(\F_p\tensor_RQ\trf_*)$.
  \begin{proof}[Proof of Claim 1]
    This is the kernel of the map
    \begin{equation*}
      \F_p\tensor_R Q\trf_*: \F_p\tensor_RQH_*(Q\Sigma\CP^\infty_+)
      \to \F_p\tensor_RQH_*(Q_0S^0)^{(-)}
    \end{equation*}
    Clearly, the natural map $JH_*(\Sigma\CP^\infty_+)\to
    k\tensor_RQH_*(Q\Sigma\CP^\infty_+)$ is an isomorphism, and by
    Lemma~\ref{lem:4.4} we get that $\F_p\tensor_RQH_*(QS^0)^{(-)}$ is
    spanned by $\{\beta Q^s\iota\mid s\geq 1\}$.  Thus Claim 1
    follows from Theorem~\ref{thm:MMM}.
  \end{proof}\noindent
  \emph{Claim 2}: $\Cok(T^\RR_1(\F_p,Q\trf_*))$ is concentrated in
  $\deg_\beta = -1$ and $\deg_\beta = -2$.
  \begin{proof}[Proof of Claim 2]
    We will compute $T^\RR_1(\F_p,Q\trf_*)$ using suitable free
    resolutions.  For brevity, write $V = JH_*(\Sigma\CP^\infty_+)$.
    By Corollary~\ref{cor:3.16} we may consider $Q\trf_*$ as a map
    from $D'V$ onto $D'\F_p^{(-)}$.  Let $W\subseteq (\DD'\F_p)^{(-)}$
    denote the subspace with basis $\{\beta Q^{s_1}Q^{s_2}\dots Q^{s_k} \mid
    s_1 \geq 1, s_2, \dots, s_k \geq 0\}$.  In the diagram
    \begin{equation}
      \label{eq:21}
      \xymatrix{
        &&W\ar[d]& D'V \ar[d]^{Q\trf_*}&\\
        0\ar[r] & {(\langle \AA\rangle\cdot \DD'\F_p)^{(-)}}\ar[r] &
        {\DD'\F_p^{(-)}}\ar[r] & D'\F_p^{(-)}\ar[r] &0
      }
    \end{equation}
    in which the lower exact sequence is an instance of
    \eqref{eq:exactseq2}, we may choose a lifting $\rho: W\to D'V$
    since $Q\trf_*$ is surjective.  Writing $V = V_0\oplus V_1$ where
    $V_0 = \textrm{span}\{a_s \mid s\equiv -1 \pmod{2(p-1)}\}$ and
    $V_1 = \textrm{span}\{a_s \mid s\not\equiv -1 \pmod{2(p-1)}\}$, we
    may choose the lifting $\rho$ to have $\rho(W)\subseteq D'V_0$
    since $D'V = D'V_0\oplus D'V_1$ and since $Q\trf_*$ vanishes on
    $D'V_1$.  We may also choose the lifting to have $\rho(\beta Q^s)
    = a_{2s(p-1)-1}$ and extend \eqref{eq:21} to the following exact
    diagram
    \begin{equation}
      \label{eq:22}
      \xymatrix{
        0\ar[r]& {\Ker(\rho)} \ar[d]^\sigma \ar[r]^j &{\DD'
        W}\ar[d]\ar[r]^\rho & D'V_0 \ar[d]^{Q\trf_*}\ar[r] &0\\
        0\ar[r] & {(\langle \AA\rangle\cdot \DD'\F_p)^{(-)}}\ar[r]^-{i} &
        {\DD'\F_p^{(-)}}\ar[r] & D'\F_p^{(-)}\ar[r] &0
      }
    \end{equation}
    Note that the middle map in \eqref{eq:22} is an isomorphism.
    
    Next we apply the functor $\F_p\tensor_\RR-$ to \eqref{eq:22}.  This
    gives a diagram involving the left derived functor $T^\RR_1(\F_p,-) =
    L_1(\F_p\tensor_\RR-)$.  This functor vanishes on the middle part of
    \eqref{eq:22} since these (isomorphic) objects are free.  Thus, a
    part of the induced diagram looks like this
    \begin{equation}\label{eq:33}
      \xymatrix{
        0\ar[r] & {T^\RR_1(\F_p,D'V_0)} \ar[d]^{T^\RR_1(\F_p,Q\trf_*)}
        \ar[r] & {\F_p\tensor_\RR \Ker\rho}\ar[d]^{\sigma_*} \ar[r]^{j_*}
        & {\F_p\tensor_\RR \DD'W} \ar[d]^{\isom}\\
        0 \ar[r] & {T^\RR_1(\F_p,D'\F_p^{(-)})} \ar[r] & {\F_p\tensor_\RR
        (\langle\AA\rangle\cdot\DD'\F_p)^{(-)}} \ar[r]^-{i_*} &
        {\F_p\tensor_\RR(\DD'\F_p)^{(-)}}
      }
    \end{equation}
    where a star in subscript is shorthand for $\F_p\tensor_\RR-$ on
    morphisms.  Thus we have represented $T^\RR_1(\F_p,D'V_0)$ and
    $T^\RR_1(\F_p,\DD'\F_p^{(-)})$ as the kernels of $j_*$ and $i_*$, and
    the map $T^\RR_1(\F_p,Q\trf_*)$ as the restriction of $\sigma_*$.
    
    To calculate the cokernel of $T^\RR_1(\F_p,Q\trf_*)$ and to prove
    Claim 2, note that
    \begin{equation*}
      (\langle\AA\rangle\cdot\DD'\F_p)^{(-)} =
      \RR^{(-)}\cdot\AA^{(0)}\cdot\DD'\F_p^{(0)} +
      \RR\cdot\AA^{(-)}\cdot\DD'\F_p^{(0)} + \RR\cdot\AA\cdot\DD'\F_p^{(-)}
    \end{equation*}
    This is generated over $\RR$ by the subspace
    \begin{equation}\label{eq:23}
      \RR^{(-1)}\cdot\AA^{(0)}\cdot\DD'\F_p^{(0)} +
      \AA^{(-)}\cdot\DD'\F_p^{(0)} + \AA\cdot\DD'\F_p^{(-)}
    \end{equation}
    The corresponding $\RR$-indecomposable classes will span
    $\F_p\tensor_\RR (\langle\AA\rangle\cdot \DD'\F_p)^{(-)}$ as a
    vectorspace, and since the first and the second term in
    \eqref{eq:23} has $\deg_\beta \in\{-1,-2\}$, it suffices to prove
    that the last term $\AA\cdot\DD'\F_p^{(-)}$ does not contribute to
    the cokernel of $T^\RR_1(\F_p,Q\trf_*)$.

    To this end, notice that $\AA\cdot\DD'\F_p^{(-)}$ corresponds to
    $\AA\cdot\DD'W$ under the middle isomorphism in~\eqref{eq:22}, and
    that $\AA\cdot\DD'W$ is in the kernel of $\rho$ since the action
    of $\AA$ is trivial in $D'V$.  Notice also that $\AA\cdot\DD'W$
    vanishes under the projection $\DD'W\to \F_p\tensor_\RR\DD'W$ and
    thus by exactness of~\eqref{eq:22} and~\eqref{eq:33} the classes
    corresponding to $\AA\cdot\DD'\F_p^{(-)}$ in $\F_p\tensor_\RR
    (\langle\AA\rangle\cdot\DD'\F_p)^{(-)}$ lifts all the way to
    $T^\RR_1(\F_p,D'V_0)$ and therefore does not contribute to the
    cokernel of $T^\RR_1(\F_p,Q\trf_*)$.
  \end{proof}
  Now Theorem~\ref{thm:4.3} follows from the exact
  sequence~\eqref{eq:20} and the Claims above.
\end{proof}

\section{Homology of $\Omega^\infty\Sigma\CPl$}
\label{sec:homol-omeg}

The spectral sequence associated to the fibration~\eqref{eq:24} has
\begin{equation}
  \label{eq:25}
  E^2 = \Cotor^{H_*(Q_0S^0)}(H_*(Q\Sigma\CP^\infty_+),\F_p) \Rightarrow H_*(\Omega^\infty\Sigma\CPl)
\end{equation}
By Proposition~\ref{prop:Cotorsplit} the $E^2$-term splits as
\begin{equation}
  \label{eq:26}
  E^2 \isom \Cotor^{H_*(Q_0S^0)\moddd\trf_*}(\F_p,\F_p) \tensor
  H_*(Q\Sigma\CP^\infty_+)\modd\trf_*
\end{equation}
In this section, $p$ is odd so after localising the
fibration~\eqref{eq:24}, the base-space is simply connected and the
spectral sequence converges.

As explained in the introduction, we will first prove
Theorem~\ref{thm-2} about the coalgebra structure on
$H_*(Q_0S^0)\moddd\trf_*$, or, equivalently the algebra structure of
$H^*(Q_0S^0)\modd\trf^*$, and then use this to prove that the spectral
sequence \eqref{eq:25} collapses.  Then a close examination of the
$E^\infty = E^2$ term will prove Theorem~\ref{thm-1}.

\subsection{The Hopf algebra cokernel of $\trf_*$}
\label{sec:hopf-algebra-cokern}

To state the results, let us introduce a bigrading of $H_*(QS^0)$.
Recall that $H_*(QS^0)$ is the free commutative algebra on the set
\begin{equation*}
  \{Q^I\iota\mid \text{$I$ admissible, $e(I)+b(I) > 0$}\}
\end{equation*}
Make it a bigraded algebra by setting $\deg_\beta(Q^I\iota) =
\deg_\beta(Q^I)$.  By the Cartan formula for the coproduct we get that
the subalgebra $H_*(QS^0)^{(0)}$ is a Hopf subalgebra, but notice that
$H_*(QS^0)$ is not a bigraded $R$-module because of the
relation~\eqref{eq:6}.  We are now ready to prove Theorem~\ref{thm-2}
in the case $p>2$.

\begin{proof}[Proof of Theorem~\ref{thm-2} for $p>2$]
  We first prove that the composition
  \begin{equation}\label{eq:27}
    H_*(QS^0)^{(0)} \to H_*(QS^0) \to H_*(QS^0)\moddd\trf_*
  \end{equation}
  is an isomorphism of Hopf algebras.

  With the bigrading introduced above, we have $H_*(QS^0) =
  H_*(QS^0)^{(0)} \oplus H_*(QS^0)^{(-)}$ where the first summand
  is a subalgebra and the second is an ideal.  Since
  $\IM(\trf_*)\subseteq \F_p\oplus H_*(QS^0)^{(-)}$, the
  composition~\eqref{eq:27} is injective.

  To see surjectivity, note that $Q(H_*(QS^0)\moddd\trf_*) =
  \Cok(Q\trf_*)$ since $Q$ is right exact.  By Theorem~\ref{thm:imF}
  we have $\IM(Q\trf_*) = QH_*(QS^0)^{(-)}$ and hence $\Cok(Q\trf_*)
  = (QH_*(QS^0))^{(0)} = Q(H_*(QS^0)^{(0)})$.
  
  To prove that $H_*(Q_0S^0)\moddd\trf_*$ is dual to a polynomial,
  notice that we have $H_*(Q_0S^0)\moddd\trf_* \isom
  H_*(Q_0S^0)^{(0)}$ and that it suffices to prove that $\lambda:
  H_*(Q_0S^0)^{(0)} \to H_*(Q_0S^0)^{(0)}$ is surjective.  $\lambda$
  is given by the dual Steenrod operations: If $\deg(x) = 2ps$,
  $\lambda x = \mathcal{P}^s_*x$.  By the Nishida relations
  (\cite[Theorem 1.1 (9)]{CLM}), one gets $\lambda Q^{ps} =
  Q^s\lambda$ and thus
  \begin{equation*}
    \lambda(Q^{ps_1}Q^{ps_2}\dots Q^{ps_k}[1]*[-p^k]) =
    Q^{s_1}Q^{s_2}\dots Q^{s_k}[1]*[-p^k]
  \end{equation*}
  Thus $\lambda$ hits the generators of $H_*(Q_0S^0)^{(0)}$ and since
  it is a map of algebras, it is surjective.
\end{proof}

\subsection{The spectral sequence}
\label{sec:spectral-sequence}
We are now ready to compute the $E^2$-term of the spectral
sequence~\eqref{eq:25} and to prove that it collapses at the
$E^2$-term.

\begin{thm}\label{thm:5.3}
  The spectral sequence collapses at the $E^2$-term.  The $E^2$-term
  is given by
  \begin{equation*}
    E^2 = H_*(Q\Sigma\CP^\infty_+)\modd\trf_* \tensor
    E[s^{-1}P(H_*(QS^0)\moddd\trf_*)]
  \end{equation*}
  as a Hopf algebra.
\end{thm}
\begin{proof}
  We need to identify the factor
  $\Cotor^{H_*(Q_0S^0)\moddd\trf_*}(\F_p,\F_p)$ in the
  splitting~\eqref{eq:26} of the $E^2$-term.  By
  Theorem~\ref{thm-2}, the dual algebra $H^*(Q_0S^0)\modd\trf^*$ is
  polynomial and hence by Corollary~\ref{lem:Cotorstruktur} we get
  \begin{equation*}
    \Cotor^{H_*(Q_0S^0)\moddd\trf_*}(\F_p,\F_p) \isom
    E[s^{-1}P(H_*(QS^0)\moddd\trf_*)]
  \end{equation*}
  as claimed.

  In this $E^2$-term, primitives and generators are concentrated in
  bidegrees $(0,*)$ and $(-1,*)$ and hence by Lemma~\ref{lem:diffHopf}
  there can be no non-zero differentials in the spectral sequence.
\end{proof}

\begin{proof}[Proof of Theorem~\ref{thm-1}]
  By Theorem~\ref{thm:5.3} we get that the sequence
  \begin{equation*}
    \xymatrix{
      {\F_p} \ar[r] & {H_*(\CPll)\modd\omega_*} \ar[r] & {H_*(\CPll)}
      \ar[r]^-{\omega_*} & {H_*(Q\Sigma\CP^\infty_+) \modd \trf_*}
      \ar[r] & {\F_p}
    }
  \end{equation*}
  is exact (i.e.\ $\omega_*$ is onto).

  To identify $H_*(\CPll)\modd\omega_*$ recall that the spectral
  sequence defines a filtration $F_0 \supseteq F_{-1} \supseteq \dots$
  on $H_*(\CPll)$ and hence on $H_*(\CPll)\modd\trf_*$ and an
  isomorphism of graded vectorspaces
  \begin{equation*}
    s^{-1}P(H_*QS^0\moddd\trf_*) \to F_{-1}(H_*(\CPll)\modd\trf_*)/F_{-2}
  \end{equation*}
  Choosing any lifting
  \begin{equation*}
    \xymatrix{
      {s^{-1}P(H_*QS^0\moddd\trf_*)} \ar[r]^-{l} &
      {H_*(\CPll)\modd\trf_*}
    }
  \end{equation*}
  we will get an isomorphism of algebras
  \begin{equation*}
    E[s^{-1}P(H_*QS^0\moddd\trf_*)] \to H_*(\CPll)\modd\trf_*
  \end{equation*}
  and since $H_*(\CPll)\modd\omega_*$ is a Hopf algebra,
  Theorem~\ref{thm:mmses} defines a unique choice of lifting $l$ into
  $P(H_*(\CPll)\modd\omega_*)$.

  The splitting follows from Lemma~\ref{lem:splitHopf} below.
\end{proof}

\begin{lem}\label{lem:splitHopf}
  Let
  \begin{equation*}
    \xymatrix{
      {\F_p}\ar[r]& A\ar[r]& B\ar[r]^\pi &C\ar[r]& {\F_p}
    }
  \end{equation*}
  be a short exact sequence of Hopf algebras.  If either $A$ or $C$ is
  exterior, the sequence is split exact in the category of Hopf algebras.
\end{lem}
\begin{proof}
  Assume $C$ is exterior.  Then by Theorem~\ref{thm:mmses} we have
  that $PC \isom QC$ and the diagram
  \begin{equation*}
    \xymatrix{
      PB\ar[r]^{P\pi} \ar[d] & PC \ar[r]\ar[d]^{\isom}&0\\
      QB \ar[r] & QC\ar[r]& 0
      }
  \end{equation*}
  is exact since $Q(-)$ is right exact.  Thus $PB\to PC$ is surjective
  and a choice of splitting $PC\to PB$ of $P\pi$ induces a splitting
  $C\isom E[PC] \to B$ of $\pi$.

  The case where $A$ is exterior follows by duality.
\end{proof}

\begin{cor}
  The vectorspace
  \begin{equation*}
    \Ker(P\omega_*) = P(H_*(\CPll)\modd\trf_*) =
    Q(H_*(\CPll)\modd\omega_*)
  \end{equation*}
  is concentrated in degrees $\equiv -1 \pmod{2(p-1)}$
\end{cor}
\begin{proof}
  This follows from Theorem~\ref{thm-2} and Theorem~\ref{thm-1}.
\end{proof}

\section{Homology of $\Omega^\infty\CPl$}
\label{sec:homol-omeg-1}
The goal of this section is to prove Theorem~\ref{thm-4} and
Corollary~\ref{cor-5}.

As mentioned in the introduction, we will consider the Eilenberg-Moore
spectral sequence of the path-loop fibration over
$\Omega^\infty\Sigma\CPl$.  From the fibration~\eqref{eq:24} one
easily gets that $\pi_1(\Omega^\infty\Sigma\CPl) = \Z$ and therefore
we have a homotopy equivalence
\begin{equation*}
  \Omega^\infty\Sigma\CPl \simeq S^1 \times \Tilde\Omega^\infty\Sigma\CPl
\end{equation*}
where $\Tilde\Omega^\infty\Sigma\CPl \to \Omega^\infty\Sigma\CPl$ is
the universal covering map.  Furthermore we have
$\Omega(\Tilde\Omega^\infty\Sigma\CPl) = \Omega^\infty_0\CPl$, the
basepoint component of $\Omega^\infty\CPl$.  Similarly
$Q\Sigma\CP^\infty_+ \simeq S^1 \times \Tilde Q\Sigma\CP^\infty_+$ and
under these splittings the map $\omega$ in the fibration~\eqref{eq:24}
restricts to a map $S^1\to S^1$ of degree 2.  
Since $p$ is odd, the effect of replacing $\CPll$ and
$Q\Sigma\CP^\infty_+$ by their universal convering spaces is to remove
a factor of $H_*(S^1) = E[\sigma]$, $\sigma = [S^1]\in H_1(S^1)$, from
each of the terms $H_*(\CPll)$ and
$H_*(Q\Sigma\CP^\infty_+)\modd\trf_*$ in Theorem~\ref{thm-1}.

The Eilenberg-Moore spectral sequence associated to the path-loop
fibration over $\Tilde \Omega^\infty\Sigma\CPl$ is
\begin{equation}
  \label{eq:15}
  E^2 = \Cotor^{H_*(\Tilde\Omega^\infty\Sigma\CPl)}(\F_p,\F_p) \Rightarrow
  H_*(\Omega_0^\infty\CPl)
\end{equation}
and by Theorem~\ref{thm-1}, the $E^2$-term splits (non-canonically) as
\begin{equation}
  \label{eq:29}
  E^2 \isom \Cotor^{H_*(\Omega^\infty\Sigma\CPl)\modd\omega_*}(\F_p,\F_p) \tensor
  \Cotor^{H_*(\Tilde Q\Sigma\CP^\infty_+)\modd\trf_*}(\F_p,\F_p)
\end{equation}
More canonically there is a short exact sequence of $\Cotor$'s, but
for the following arguments we will assume that a splitting has been
chosen.

I claim it must collapse.  As before, we consider a possibly nonzero
differential $dx = y\neq 0$ with $\deg(x)$ minimal.  We will reach a
contradiction in a number of steps.  The argument is based on
Theorem~\ref{thm:4.3} and a careful analysis of degrees modulo
$2(p-1)$ in the spectral sequence.

By Theorem~\ref{thm-1} and Lemma~\ref{lem:Cotorstruktur}, the first
factor in \eqref{eq:29} is a polynomial algebra on generators of total
degree $\equiv -2 \pmod{2(p-1)}$.  To gain information about the
second factor, we map the spectral sequence~\eqref{eq:15} into the
spectral sequence of the path-loop fibration over $\Tilde
Q\Sigma\CP^\infty_+$ via the map $\omega:
\Tilde\Omega^\infty\Sigma\CPl\to \Tilde Q\Sigma\CP^\infty_+$.  This is
a map $E^r(\omega)$ of spectral sequences whose restriction to the
first factor in the splitting~\eqref{eq:29} is zero, and whose
restriction to the second factor in~\eqref{eq:29} is induced by the
inclusion $H_*(\Tilde Q\Sigma\CP^\infty_+)\modd\trf_* \to H_*(\Tilde
Q\Sigma\CP^\infty_+)$.  The next lemma says that this second factor
in~\eqref{eq:29} injects under $E^2(\omega)$.

\begin{lem}\label{lem:6.1}
  Let $f: A\to B$ be an injection of primitively generated Hopf
  algebras.  Then $\Cotor^f(\F_p,\F_p): \Cotor^A(\F_p,\F_p)\to
  \Cotor^B(\F_p,\F_p)$ is also injective.
\end{lem}
\begin{proof}
  By Theorem~\ref{thm:mmses}, $A^*$ and $B^*$ are tensor products of
  exterior algebras and polynomial algebras truncated at height $p$.
  Thus we can split $f^*: B^*\to A^*$ in the category of algebras
  (since a splitting can be chosen on the generators of $A^*$).
  Dually, $f: A\to B$ is split injective as a map of coalgebras and
  thus $\Cotor^f(\F_p,\F_p)$ is injective.
\end{proof}
\begin{cor}
  \label{cor:6.2}
  Relative to the splitting \eqref{eq:29}, a differential $dx=y\neq 0$
  with $x$ of minimal degree will have $x$ in the right factor and $y$
  in the left.
\end{cor}
\begin{proof}
  Recall that $P$ and $Q$ are additive: $P(A\tensor B) = PA \oplus
  PB$ and $Q(A\tensor B) = QA\oplus QB$.  Thus $x$ and $y$ does not
  contain products between the two factors in \eqref{eq:29}.
  
  Since $y$ is primitive and in bidegree $(\leq -3,*)$, it must be of
  even total degree by Corollary~\ref{lem:Cotorstruktur}, and thus $x$
  is of odd total degree.  By Theorem~\ref{thm-1} this is only
  possible if $x$ is in the right factor.

  By Lemma~\ref{lem:6.1}, the right factor injects into the spectral
  sequence of $Q\Sigma\CP^\infty_+$, and since  all differentials
  vanish in this spectral sequence, $y$ must map to 0 there, and hence
  $y$ is in the left factor.
\end{proof}

The remaining part of the collapse proof is to eliminate the
possibility of differentials from the right factor to the left.  This
is the hardest part of the proof, the main ingredient of which is
Theorem~\ref{thm:4.3}.
\begin{thm}
  \label{thm:6.4}
  The spectral sequence \eqref{eq:15} collapses.
\end{thm}
\begin{proof}
  Assume there is a differential $dx=y\neq 0$ with $\deg(x)$ minimal.
  Then $y$ is a primitive element in
  $\Cotor^{H_*(\Tilde\Omega^\infty\Sigma\CPl)}(\F_p,\F_p)$.  By
  Corollary~\ref{cor:6.2} and Corollary~\ref{lem:Cotorstruktur} it is
  of the form
  \begin{equation*}
    y = (s^{-1}z)^{p^k}
  \end{equation*}
  for a $z \in P(H_*(\Tilde\Omega^\infty\Sigma\CPl)\modd\omega_*)$.
  By Theorem~\ref{thm-1} we must have $\deg(z) \equiv
  -1\pmod{2(p-1)}$.  Write
  \begin{equation*}
    \deg(z) = 2n(p-1)-1
  \end{equation*}
  Then
  \begin{equation*}
    \deg y = p^k(2n(p-1)-2) = 2p^k(n(p-1)-1) \equiv -2 \pmod{2(p-1)}
  \end{equation*}
  and thus $\deg x \equiv -1 \pmod{2(p-1)}$ because the differential
  has degree $-1$.  By Proposition~\ref{lem:diffHopf} we get that $x$
  corresponds to a minimal element in the cokernel of $\sigma_*:
  QH_*(\Omega^\infty_0\CPl) \to PH_*(\Tilde\Omega^\infty\Sigma\CPl)$,
  of degree $\equiv 0 \pmod{2(p-1)}$.  By Corollary~\ref{cor:6.2}, $x$
  is also a minimal element in the cokernel of the composition
  \begin{equation*}
    \xymatrix{
      {QH_*(\Omega^\infty_0\CPl)} \ar[r]^{\sigma_*} &
      {PH_*(\Tilde\Omega^\infty\Sigma\CPl)} \ar[r]^-{P\omega_*} &
      {P(H_*(\Tilde Q\Sigma\CP^\infty_+)\modd\trf_*)}
    }
  \end{equation*}
  By minimality this element is not a $p$th power and hence maps to a
  non-zero element of $Q(H_*(\Tilde Q\Sigma\CP^\infty_+)\modd\trf_*)$.
  Again by minimality, and because the loop suspension $\sigma_*$ is
  $R$-linear, this element is $R$-indecomposable and hence since
  $\sigma_*$ has degree 1, $x$ will map to a nonzero element of degree
  $\equiv 0\pmod{2(p-1)}$ in
  \begin{equation*}
    \F_p\tensor_R Q(H_*(\Tilde Q\Sigma\CP^\infty_+)\modd\trf_*)
  \end{equation*}
  in contradiction with Theorem 4.3.
\end{proof}
Theorem~\ref{thm-4} no follows from Theorem~\ref{thm:6.4} and
Proposition~\ref{cor:mindiff}.

\begin{proof}[Proof of Corollary~\ref{cor-5}]
  This is completely analogous to the inductive step in the classical
  calculations of homology of $QX$ or of cohomology of $K(\F_p,n)$.
  We sketch the details.

  Consider the Leray-Serre spectral sequence 
  \begin{equation}
    \label{eq:4}
    E^2 = H_*(\CPll)\tensor H_*(\CPlll) \Rightarrow
    H_*(\mathit{point})
  \end{equation}
  Since $\sigma_*$ is onto, we can pick a basis $B\subseteq
  PH_*(\CPll)$ and for each $x\in B$ pick an element $\tau x \in
  H_*(\CPlll)$ with $\sigma_*(\tau x) = x$.  We can now form a model
  spectral sequence 
  \begin{equation*}
    \Tilde E^2 = \bigotimes_{x\in B} E^r(x)
  \end{equation*}
  where, if $x\in B$ has odd degree,
  \begin{equation*}
    E^2(x) = E[x] \tensor \F_p[\tau x]
  \end{equation*}
   with the differential determined by requiring that $x$ transgresses
   to $\tau x$.  If $x$ has even degree $\deg(x) = 2s$, we set
   \begin{equation*}
     E^2(x) = \F_p\{1,x,\dots, x^{p-1}\} \tensor E[\tau x] \tensor
     \F_p[\beta Q^s(\tau x)]
   \end{equation*}
   with the differential determined by requiring that $x$ transgresses
   to $\tau x$ and that $x^{p-1} \tensor \tau x$ transgresses to
   $\beta Q^s (\tau x)$.

   The choices of $\tau x\in H_*(\CPll)$ determines a map of spectral
   sequences $\Tilde E^r \to E^r$
   and the comparison theorem implies that it is an isomorphism and
   then Corollary~\ref{cor-5} follows.
\end{proof}
\section{The case $p=2$}\label{sec:case-p=2}

At the prime 2, the calculation of $H_*(\Omega^\infty\CPl)$ and
$H_*(\Omega^\infty\Sigma\CPl)$ can also be made.  Some details are
quite different however.  In particular, we will use the looped
fibration
\begin{equation}
  \label{eq:34}
  \Omega^\infty\CPl\to Q(\CP^\infty_+) \to \Omega QS^0
\end{equation}
to compute $H_*(\Omega^\infty\CPl)$, instead of the path-loop
fibration over $\Omega^\infty\Sigma\CPl$.  At $p=2$ our base spaces in
the fibrations are no longer simply connected.  The following lemma
deals with this
\begin{lem}\label{lem:splitspace}
  As spaces we have
  \begin{align*}
    QS^0 & \simeq \Z \times \RP^\infty \times \tilde Q_0S^0\\
    \Omega QS^0 &\simeq \Z/2 \times \RP^\infty\times \Tilde\Omega_0 QS^0
  \end{align*}
  where $\Tilde X\to X$ denotes the universal covering.
\end{lem}
\begin{proof}
  Let $X$ be an $(n-1)$-connected $H$-space with $\pi_n(X) = G$. There
  is an $H$-map $X\to K(G,n)$ inducing an isomorphism in $\pi_n$ and
  with fibre the $n$-connected cover $X\langle n\rangle$. If one can
  find a map $K(G,n)\to X$ inducing an isomorphism in $\pi_n$, this
  map will give a splitting $X \simeq X\langle n\rangle\times K(G,n)$.
  
  For $n=0$ this is automatic.
  
  For $X=Q_2S^0\simeq Q_0S^0$, $\pi_1(X) = \Z/2$ and the definition
  of the Dyer-Lashof operation $Q^1\iota\in H_1(Q_2S^0;\F_2)$ gives a map
  \begin{displaymath}
    \RP^\infty = B\Z/2\to Q_2S^0 \simeq Q_0S^0
  \end{displaymath}
  inducing an isomorphism in $H_1$ and thus by the Hurewicz theorem an
  isomorphism in $\pi_1$ and the splitting of $QS^0$ follows.
  
  For $X=\Omega_0Q_0S^0$, $\pi_1(X) = \Z/2$. The Hopf map gives an
  infinite loop map $\eta: Q(S^1)\to Q_0S^0$. I claim it is nonzero in
  $\pi_2$. To see this it suffices to show that
  $(\eta\langle1\rangle)_*$ is nonzero in $H_2$ which can be seen as
  follows. Let $\sigma\in H_1(QS^1)$ be the fundamental class. Since
  $QS^1 \simeq S^1\times QS^1\langle1\rangle$, the element
  $Q^1\sigma\in H_2(QS^1)$ must be in the image from
  $H_*(QS^1\langle1\rangle)$. Since $\eta_*(Q^1\sigma) =
  Q^1(Q^1[1]*[-2]) \neq 0$, $\eta\langle1\rangle_*$ is indeed nonzero
  in $H_2$.

  Hence, $\Omega_0\eta: Q_0S^0 \to \Omega_0Q_0S^0$ is nonzero in
  $\pi_1$ and thus the composition
  \begin{displaymath}
    \RP^\infty\to Q_0S^0\to \Omega_0Q^0S^0
  \end{displaymath}
  is nonzero in $\pi_1$ and the splitting of $\Omega QS^0$ follows.
\end{proof}

Lemma~\ref{lem:splitspace} ensures that our spectral sequences has
trivial local coefficients and hence that the spectral sequences
converges.

\subsection{Recollections}
\label{sec:recollections-1}

The Dyer-Lashof algebra is slightly different at $p=2$. Let $\RR$ be
the free non-commutative algebra on the set $\{Q^s\mid s\geq 0\}$ with
$\deg(Q^s) = s$.  The Adem relation $\AA^{(0,r,0,s)}$ in
Definition~\ref{defn:Adem} still makes sense, and we let
$\AA\subseteq\RR$ be the span of the $\AA^{(0,r,0,s)}$.  The
unstability relations at $p=2$ are
\begin{equation*}
  Q^s x = 
  \begin{cases}
    x^2 & \text{if $\deg x = s$}\\
    0 & \text{if $\deg x  > s$}
  \end{cases}
\end{equation*}
and the algebra $R$ is defined from these data as before.
Corresponding to $I = (s_1,s_2, \dots, s_k)$ there is an iterated
operation $Q^I = Q^{s_1} \dots Q^{s_k}$, and this operation is called
admissible if $s_i \leq 2s_i$ for all $i$.  The definition of excess
at $p=2$ is
\begin{equation*}
  e(I) = s_1 - \sum_{j=2}^k s_j
\end{equation*}
Given
a basis $B\subseteq JH_*(X)$, then $H_*(QX)$ is the
polynomial algebra on the set
\begin{equation*}
  \mathbf{T} = \{Q^I x \mid \text{$x\in B$, $I$ admissible, $e(I) >
  \deg(x)$}\}
\end{equation*}
and similarly for $H_*(Q_0X)$.

One pleasant feature of $p=2$ is the following
\begin{lem}
  The cohomology algebra $H^*(Q_0X)$ is polynomial if the Frobenius
  $\xi: H^*(X) \to H^*(X)$ is injective.
\end{lem}
\begin{proof}
  This is because the Nishida relation $\lambda Q^{2s} = Q^s\lambda$
  makes $\lambda: H_*(Q_0X)\to H_*(Q_0X)$ surjective if $\lambda:H_*(X)
  \to H_*(X)$ is surjective.
\end{proof}
In particular, $H^*(Q_0S^0)$ and $H^*(Q_0\CP^\infty_+)$ are both
polynomial.

The calculation in Theorem~\ref{lem:Torstruktur} is valid with the
remark that $\F_2[x]/(x^2)$ must be interpreted as $E[x]$ and thus it
does not produce generators of $\Cotor$ in bidegree $(-2,*)$.  Only
truncations at height $p^n, n\geq 2$ does that.

An important difference is that for odd primes, $\Cotor^A(\F_p,\F_p)$ is
automatically a free algebra.  This is no longer true for $p=2$, since
$\Tor^{\F_2[x]}(\F_2,\F_2) = E[s^{-1}x]$, and exterior algebras are not free
in characteristic 2.
  
One consequence of the above remarks is the following
\begin{prop}
  \label{prop:suspisom}
  Let $X$ be a simply connected space with $H^*(X)$ polynomial.  Then
  $H_*(\Omega X)$ is an exterior algebra and the suspension
  \begin{equation*}
    \sigma_*: QH_*(\Omega X) \to PH_*(X)
  \end{equation*}
  is an isomorphism.  The spectral sequence
  \begin{equation*}
    \Cotor^{H_*(X)}(\F_2,\F_2) \Rightarrow H_*(\Omega X)
  \end{equation*}
  collapses.
\end{prop}
\begin{proof}
  This is because
  \begin{equation*}
    \Cotor^{H_*(X)}(\F_2,\F_2) \isom E[s^{-1}PH_*(X)]
  \end{equation*}
  has generators and primitives in bidegrees $(-1,*)$.  Together with
  Lemma~\ref{lem:diffHopf}, this proves the collapse claim and that
  $\sigma_*$ is an isomorphism.  Dually we have that
  \begin{equation*}
    \sigma^*: QH^*(X) \to PH^*(\Omega X)
  \end{equation*}
  is an isomorphism and since the image generates $H^*(\Omega X)$ as
  an algebra, $H^*(\Omega X)$ is primitively generated.  By
  Theorem~\ref{thm:mmses} we get that $H_*(\Omega X)$ is exterior.
\end{proof}
In particular this applies to $X = \Tilde Q_0S^0$.

Similarly, we have
\begin{prop}
  For any space $X$, the spectral sequence
  \begin{equation*}
    \Cotor^{H_*(\Tilde Q\Sigma X)}(k,k) \Rightarrow H_*(Q_0X)
  \end{equation*}
  collapses and the suspension
  \begin{equation*}
    \sigma_*: QH_*(QX) \to PH_*(Q\Sigma X)
  \end{equation*}
  is an isomorphism.
\end{prop}
\begin{proof}
  $\sigma_*$ is surjective since it hits $JH_*(\Sigma X)$ and since it
  is $R$-linear.  Thus by Corollary~\ref{cor:mindiff}, the spectral
  sequence must collapse.  Now $H_*(Q\Sigma X)$ is primitively
  generated, so by Theorem~\ref{thm:mmses} we get that $H^*(Q\Sigma
  X)$ is exterior and hence the spectral sequence has
  \begin{equation*}
    E^2 = \Cotor^{H_*(\Tilde Q\Sigma X)}(\F_2,\F_2) \isom
    \F_2[s^{-1}PH_*(\Tilde Q\Sigma X)]
  \end{equation*}
  Since this is free as an algebra, there are no extension problems in
  homology, and since $QH_*(Q_0X)$ is in linear bijection with
  $E^\infty_{-1,*}$, we get that $\sigma_*$ is injective.
\end{proof}

\subsection{Homology of $\Omega^\infty\CPl$}
\label{sec:homology-cpll}

The lemmas in subsection~\ref{sec:recollections-1} imply the
following diagram
\begin{equation}
  \label{eq:3}
  \begin{array}{c}
    \xymatrix{
      {QH_*(Q_0\CP^\infty_+)} \ar[r]^-{Q(\Omega_0\trf)_*} \ar[d]_{\isom} &
      {QH_*(\Omega_0 QS^0)} \ar[d]^{\isom}\\
      {PH_*(\Tilde Q\Sigma\CP^\infty_+)} \ar[r]^-{P\trf_*} &
      {PH_*(\Tilde Q_0S^0)}
    }
  \end{array}
\end{equation}
in which the vertical isomorphisms are the suspensions and in which
$H_*(\Omega_0 QS^0)$ is an exterior algebra, dual to a polynomial
algebra.

The formula for $\trf_*$ has an extra term because of the Hopf map
$\eta$.  We quote the result from \cite[Theorem 4.4]{MMM}:
\begin{thm}[\cite{MMM}]\label{thm:MMM2}
  Let $a_s\in H_*(\CP^\infty_+)$ be the generator, $s$ odd. Then
  \begin{equation*}
    Q(\trf_*)(a_s) = Q^{2s+1}\iota + Q^{s+1}Q^s \iota = Q^{2s+1}\iota
    + Q^{2s}Q^1\iota
  \end{equation*}
  \qed
\end{thm}

We shall need a lemma analogous to Lemma~\ref{lem:4.4}.
\begin{lem}\label{lem:lemma1}
  The left ideal in $R$ generated by $\{Q^{2s+1} \mid s\geq 0\}$ is
  also a right ideal.
\end{lem}
\begin{proof}
  This is completely analogous to the proof of Lemma~\ref{lem:4.4}.
  One uses the Adem relation
  \begin{equation*}
    Q^{2s}Q^{r-s} = Q^r Q^s + \sum_{i>s} \lambda_i Q^{r+s-i}Q^i
  \end{equation*}
  valid for $r \leq 2s$, for $r$ odd and $s$ even.
\end{proof}
\begin{lem}\label{lem:7.5}
  Let $b_{2s+1} \in PH_*(QS^0) = PH_*(Q_0S^0)$ be the unique primitive
  element with $b_{2s+1} - Q^{2s+1}\iota$ decomposable.  Then
  $PH_*(QS^0)$ is generated over $R$ by the set $\{b_{2s+1}\mid
  s\geq0\}$.
\end{lem}
\begin{proof}
  Let $\lambda: QH_*(Q_0S^0)\to QH_*(Q_0S^0)$ be the dual of the
  squaring.  By the Nishida relation $\lambda Q^{2s} = Q^{s}\lambda$,
  the coimage of $\lambda$ has basis
  \begin{equation*}
    \{Q^I\iota \vert \text{$I$ admissible, $e(I) > 0$,
    $2|I$}\}
  \end{equation*}
  where $2|I$ means that all entries of $I$ are even.  Thus
  Theorem~\ref{thm:mmses} implies that the image of $PH_*(QS^0) \to
  QH_*(QS^0)$ has basis
  \begin{equation*}
    \{Q^I\iota \mid \text{$I$ admissible, $e(I) > 0$,
    $2\not|I$}\}
  \end{equation*}
  and by Lemma~\ref{lem:lemma1}, this is generated over $R$ by the
  subset
  \begin{equation*}
    \{Q^{2s+1}\iota \mid s\geq 0\}
  \end{equation*}
  
  Thus the subspace of $PH_*(QS^0)$ generated over $R$ by
  $\{b_{2s+1}\mid s\geq 0\}$ contains all indecomposable primitives.
  But this generated subspace is clearly preserved by the Frobenius
  map $\xi: x\mapsto x^2$, so it contains all primitives and the claim
  follows from Theorem~\ref{thm:mmses}.
\end{proof}

We are now ready to prove the mod 2 analogue of Theorem~\ref{thm-2}.
The result is much simpler, and the extra term in
Theorem~\ref{thm:MMM2} does not give much trouble.
\begin{thm}\label{thm:Psurj}
  The map
  \begin{equation*}
    P\trf_*: PH_*(Q\Sigma\CP^\infty_+) \to PH_*(QS^0)
  \end{equation*}
  is surjective.
\end{thm}
\begin{proof}
  By the previous lemma, it suffices to prove that $Q\trf_*$ hits the
  classes $Q^{2s+1}\iota$.  Indeed, any indecomposable class
  mapping to $Q^{2s+1}\iota$ is odd-dimensional and thus
  by Theorem~\ref{thm:mmses} has a unique primitive representative that will
  map to $b_{2s+1}$.

  For $s=0$, this is immediate, since
  $\trf_*(a_1) = (Q^1\iota)*\iota^{-2}$.  For general $s$ we use the Adem
  relation $Q^{2s}Q^1 = Q^{s+1}Q^{s}$ to get
  \begin{align*}
    Q(\trf_*)(a_{2s+1}) &= Q^{2s+1}\iota + Q^{s+1}Q^s\iota\\
    &= Q^{2s+1}\iota + Q^{2s}Q^1\iota\\
  \end{align*}
  Thus we have
  \begin{equation*}
    Q(\trf_*)(a_s - Q^{2s}a_1) = Q^{2s+1}\iota
  \end{equation*}
\end{proof}
\begin{rem}
  The claim of \cite[Cor.\ 7.5]{MMM} that $\trf_*$ and thus
  $P(\trf_*)$ is \emph{injective} is incorrect. The $Q^IQ^{2r+1}$ of
  \cite[Cor.\ 7.4]{MMM} is not necessarily admissible, and in fact an
  application of the Adem relations shows that
  \begin{displaymath}
    \trf_*(Q^3 a_1 - Q^2Q^1a_1) = 0
  \end{displaymath}
\end{rem}

Together with the diagram~\eqref{eq:3}, Theorem~\ref{thm:Psurj} makes
the spectral sequence
\begin{equation}
  \label{eq:17}
  \Cotor^{H_*(\Tilde \Omega QS^0)}(H_*(\Tilde Q\CP^\infty_+),k) \Rightarrow
  H_*(\Omega_0^\infty\CPl)
\end{equation}
very simple.  We can now prove Theorem~\ref{thm-3}.
\begin{proof}[Proof of Theorem~\ref{thm-3}]
  It follows from diagram~\eqref{eq:3} and Theorem~\ref{thm:Psurj}
  that the map $Q(\Omega_0\trf_*)$ is surjective.  Therefore the
  $E^2$-term of the Eilenberg-Moore spectral sequence is
  \begin{equation*}
    E^2 = \Cotor^{H_*(\Omega_0QS^0)}(H_*(Q_0\CP^\infty_+),\F_2)
    \isom H_*(Q_0\CP^\infty_+) \modd \Omega\omega_*
  \end{equation*}
  and is concentrated on the line $E^2_{0,*}$.  Therefore it collapses
  and we get the short exact sequence
  \begin{equation*}
    \xymatrix{
      {\F_2} \ar[r] & {H_*(\CPlll)} \ar[r]^-{\Omega\omega_*} &
      {H_*(Q\CP^\infty_+)} \ar[r]^-{\Omega\partial_*} & {H_*(\Omega
        QS^0)} \ar[r] & {\F_2}
    }
  \end{equation*}
\end{proof}

\subsection{Homology of $\Omega^\infty\Sigma\CPl$}

This part of the calculation is similar to the odd primary case.  We
consider again the spectral sequence~\eqref{eq:25} with the
splitting~\eqref{eq:26}.  Notice that the fibration~\eqref{eq:24}
splits off the fibration $S^1 \to S^1 \to \RP^\infty$ and hence it has
trivial local coefficients.  As for odd primes, we need to determine
the coalgebra structure on $H_*(QS^0)\moddd\trf_*$,so we first prove
Theorem~\ref{thm-2} in the case $p=2$.
\begin{proof}[Proof of Theorem~\ref{thm-2}, $p=2$]
  Since $Q$ is right exact we have $Q(H_*(QS^0)\moddd\trf_*) =
  \Cok(Q\trf_*)$, and from the calculation in the proof of
  Theorem~\ref{thm:Psurj} it follows that the image of $Q\trf_*$
  contains all $Q^I\iota$ where $I$ has at least one \emph{odd} entry,
  and therefore that the composition in Theorem~\ref{thm-2} is
  surjective.

  To prove injectivity, consider again the dual squaring $\lambda:
  H_*(Q_0S^0) \to H_*(Q_0S^0)$.  It is a map of Hopf algebras, and
  since $\lambda Q^{2s} = Q^s\lambda$ and $\lambda Q^{2s+1} = 0$ we
  get that
  \begin{equation*}
    \lambda: H_*(QS^0)^{(0)} \to H_*(QS^0)
  \end{equation*}
  is an isomorphism.  Hence
  \begin{equation*}
    H_*(QS^0) = H_*(QS^0)^{(0)} \oplus \Ker(\lambda)
  \end{equation*}
  where the first summand is a subalgebra and the second is an ideal.
  Now the injectivity of the map in the theorem follows from the fact
  that $\Ker(\lambda)$ is an ideal and that $\IM(\trf_*)\subseteq
  \F_2\oplus \Ker(\lambda)$.
\end{proof}

\begin{thm}\label{thm:8.2}
  $H_*(Q_0S^0)\moddd\trf_*$ is dual to a polynomial algebra.
\end{thm}
\begin{proof}
  This follows since $\lambda: H_*(Q_0S^0) \to H_*(Q_0S^0)$ is
  surjective.
\end{proof}
Notice that $H^*(Q_0S^0)$ itself is polynomial.  This is in contrast
to the odd primary case, where only the subalgebra
$H^*(Q_0S^0)\modd\trf^*\subseteq H^*(Q_0S^0)$ is polynomial.

\begin{proof}[Proof of Theorem~\ref{thm-3}, $p=2$]
  Completely as for odd primes, Theorem~\ref{thm:8.2} makes the
  spectral sequence collapse, and the collapse gives a short exact
  sequence of Hopf algebras
  \begin{equation}\label{eq:5}
    \xymatrix{
      {\F_2} \ar[r] & {H_*(\CPll)\modd \omega_*} \ar[r]& {H_*(\CPll)}
      \ar[r]^-{\omega_*} & {H_*(Q\Sigma\CP^\infty_+) \modd\trf_*}
      \ar[r] &{\F_2}
    }
  \end{equation}
  Since $H_*(Q\Sigma\CP^\infty_+)$ is primitively generated, so is
  $H_*(Q\Sigma\CP^\infty_+) \modd \trf_*$.  Hence the sequence is
  split if and only if $P(\omega_*)$ is surjective.  We have the
  diagram
  \begin{equation*}
    \xymatrix{
      {QH_*(\CPlll)} \ar[r] \ar[d]^-{\sigma_*} &
      {QH_*(Q\CP^\infty_+)} \ar[r]\ar[d]^{\isom} &
      {QH_*(\Omega QS^0)} \ar[r]\ar[d]^{\isom} & 0\\
      {PH_*(\CPll)} \ar[r]  &
      {PH_*(Q\Sigma\CP^\infty_+)} \ar[r]^-{P\trf_*} &
      {PH_*(QS^0)} \ar[r] & 0\\
    }
  \end{equation*}
  which we know is exact except possibly at
  $PH_*(Q\Sigma\CP^\infty_+)$.  But it follows from the rest of the
  diagram that it is also exact at $PH_*(Q\Sigma\CP^\infty_+)$.  Since
  $\Ker(P\trf_*) = P(H_*(Q\Sigma\CP^\infty_+)\modd\trf_*)$ we get that
  the sequence~\eqref{eq:5} splits.
\end{proof}

\hspace{10mm}

\end{document}